\documentclass[11pt]{elsarticle}

\usepackage{amssymb}

\usepackage{amsthm}
\usepackage{graphics}
\usepackage{caption2}
\usepackage{psfrag}
\usepackage{showlabels}
\usepackage{color}
\usepackage{amsmath}
\usepackage{bm}
\usepackage{mathtools}
\usepackage{tensor}
\usepackage{subfigure}
\usepackage[normalem]{ulem}
\usepackage[linesnumbered,boxed,ruled,commentsnumbered]{algorithm2e}
\usepackage{graphicx}
\usepackage{epstopdf}
\usepackage{array}

\usepackage{amsfonts}
\usepackage{amssymb}
\usepackage{graphicx}
\usepackage{multirow}
\usepackage{exscale}
\usepackage{relsize}
\usepackage[normalem]{ulem}
\usepackage{color}
\usepackage{fullpage}
\usepackage{lineno}

\newcommand{\ds}{\displaystyle}
\newcommand{\f}{\frac}

\newcommand{\x}{\mathbf{x}}

\newcommand{\be}{\begin{equation}}
\newcommand{\ee}{\end{equation}}
\newcommand{\ba}{\begin{array}}
\newcommand{\ea}{\end{array}}

\newtheorem{remark}{Remark}

\journal{Computer Methods in Applied and Mechanical Engineering}

\begin{document}
	
\begin{frontmatter}

\title{A fast computational framework for the linear bond-based peridynamic model}

\author{Chenguang Liu$^a$ \quad  Hao Tian$^a$ \quad  Wai Sun Don$^a$ \quad Hong Wang$^b$}

\address{$^a$ School of Mathematical Science, Ocean University of China, Qingdao, Shandong 266100, China \\ $^b$Department of Mathematics, University of South Carolina, Columbia, South Carolina 29208, USA}

\begin{abstract}
Peridynamic (PD) theory is significant and promising in engineering and materials science; however, it imposes challenges owing to the enormous computational cost caused by its nonlocality.
Our main contribution, which overcomes the restrictions of the existing fast method, is a general computational framework for the linear bond-based peridynamic models based on the meshfree method, called the matrix-structure-based fast method (MSBFM), which is suitable for the general case, including 2D/3D problems, and static/dynamic issues, as well as problems with general boundary conditions, in particular, problems with crack propagation.
Consequently, we provide a general calculation flow chart. The proposed computational framework is practical and easily embedded into the existing computational algorithm.
With this framework, the computational cost is reduced from $O(N^2)$  to $O(N\log N)$, and the storage request is reduced from $O(N^2)$ to $O(N)$, where N is the degree of freedom. Finally, the vast reduction of the computational and memory requirement is verified by numerical examples.

\end{abstract}


\end{frontmatter}

\section{Introduction}

Classical continuum mechanics are expressed as a partial differential equation, which is challenging to describe models with discontinuities. As a result, peridynamics(PD), as proposed by Silling\cite{SA}, is an integral-type nonlocal model and can provide a general theory for solving problems in the form of discontinuities.
Over the past few decades, the effectiveness of PD has attracted extensive research conducted on modeling methods, numeral techniques, and applications. In this paper, we focus on the bond-based peridyanmics \cite{AXS}, which is an early version of peridynamics and can be applied to isotropic materials, with Poisson's ratio of 1/4 for plane strain and 1/3 for plane stress. Later the ordinary state-based and non-ordinary state-based PD were proposed to eliminate the constraints of fixed Poisson's ratios on materials\cite{SA2}. The PD models have been frequently used in many practical problems. A range of cutting-edge applications can be found in composite material deformation\cite{GRUA,YNE,VB}, corrosion\cite{SZF,ZSJ,JSM,SZS,SZJF}, damage prediction\cite{CSE,XXXE2,DW,MKA,PZ} and crack simulation for a variety of materials\cite{WNS,YXG,SE,SYL,PS}.

There have been much research directed at developing numerical methods that can solve the PD models, including meshfree methods, finite difference, finite element methods, and collocation methods\cite{XQ,TCP,XM,SilAsk,JY,JMY,SY}.Among other research in the literature, explored in\cite{XQ,XQ2,QJ}, the asymptotically compatible schemes retained a limit behavior that makes the limit of the zero-horizon of the nonlocal operator become the local differential operator, providing a consistency between local and non-local models. However, these methods are very restrictive due to the vast computational cost caused by nonlocality of PD. The increasing computation cost limits the application of PD theory, especially for multidimensional cases. Lots of efforts have been made to overcome this issue. A class of coupled method was introduced to accelerate the PD simulation, which utilizes PD on the area around the cracks and classical mechanics on the rest area \cite{MT,TM,JFZJ, FGY,HHHY,MPRM}.  A fast method based on the convolution structure of PD models\cite{SF,SL} is proposed to accelerate the simulation. A super-fast peridynamic model\cite{DBVK} based on decreasing the number of inner loop operations is also introduced to overcome this difficulty.  
The work above also gives us some inspiration for this article.

In the literature, a class of fast methods utilizing the structure of stiff matrix are booming, which can reduce the computational cost from  $O(N^2)$ to $O(N\log N)$ without loss of accuracy. In 2010, a fast method\cite{HH} based on the Toeplitz structure of stiff matrices was proposed, hence solving the 1D static linear bond-based peridynamics. Subsequently, a fast collocation method based on TBT matrix structure was given for 2D nonlocal diffusion models in reference\cite{HH2}, which can be thought of a approximation model of scalar-valued. A fast collocation method for the 2D static linear bond-based peridynamics with volume boundary conditions was investigated in \cite{XH} , where we use an equivalent but more effective way to evaluate. In 2017, a fast method was also presented to solve nonlocal diffusion models with variable coefficients\cite{CH}, and a discontinuous method was discussed to solve linear bond-based PD models with discontinuous solves\cite{HAH}. In 2020, a fast algorithm with preconditioned processing was proposed\cite{XXAH,XAH}, accelerating the convergence of the iterative method.  Although the research above has significantly contribute to applying fast matrix structure based methods in the simulation of PD, there are still several pending issues. The first one is that all this research is proposed for 1D or 2D static models. The second one is that there need to be methods taking into account volume constrained boundary conditions. The third one is that all these research is developed for models with no cracks. A natural question comes into the work here. To fill this gap, we give a fast matrix-based method(MSBFM), which is the main contribution of this article.

In this paper, we offer an insightful look at a very general setting. We propose a fast matrix-based method(MSBFM) for solving 2D/3D linear bond-based peridynamic models. By establishing the relationship between the stiffness matrix and the Toeplitz-Block-Toeplitz(TBT) matrix, we focos on the matrix decomposition, stiff matrix, one for TBT matrix and another for sparse matrix. This reduces the structural limitations of the matrices and hence a more efficient method applied to general boundary and cracks, and dramatically reduces the amount of computation and storage from $O(N^2)$ to $O(N\log N)$ and from $O(N^2)$ to $O(N)$, respectively. Meanwhile, the results are for models in 2D as well as in 3D. Numerical experiments verify the accuracy of the method.

The following articles are organized as follows. In section 2, we reviewed the linear model of the peridynamic and meshfree methods. In section 3, we analyze the matrix structure of the two-dimensional problem and give an accelerated process based on the matrix structure. In section 4, the matrix structure of the 3D model is discussed, and the MSBFM method is introduced. In section 5, the accuracy and acceleration effect of the MSBFM method are shown by numerical examples.

\section{Fundamentals of linear bond-based peridynamics and its discretization}
The bond-based peridynamics, is a reformulation for classical continuum solid mechanics by the integral form instead of partial differential equations. Typically, with sufficiently small displacement, bond-based peridynamics can be approximated as linear bond-based peridynamics. In this section, we mainly review the linearized version of the  bond-based peridynamic mode, addressed in this study, and the discretization to solve the model. 
\subsection{Linear bond-based peridynamics}
The equation of motion for the linear bond-based peridynamics with prescribed volume boundary condition can be defined as follows:
 \begin{equation}\label{mat:e1}
 \begin{cases}
\rho \ddot{\mathbf{u}}(\mathbf{x}, t)= \ds \int_{\mathbf{\mathcal{B}_{\delta}(\mathbf{x})}\cap\Omega} \mu\mathbf{C}(\mathbf{\mathbf{x^{'}}}-\mathbf{x})(\mathbf{u}\left(\mathbf{x}^{\prime}, t\right)-\mathbf{u}(\mathbf{x}, t)) dV_{\mathbf{\mathbf{x}^{\prime}}}+\mathbf{b}(\mathbf{x},t)&\text{$\mathbf{x}\in \Omega_s$},\\
\mathbf{u}(\mathbf{x},t)=\mathbf{h}(\mathbf{x,t})&\text{$\mathbf{x}\in \Omega_c$}.
\end{cases}
\end{equation}
where $\rho$ is the mass density, $\Omega=\Omega_s\cup\Omega_c$ is the spatial domain, $\mathcal{B_{\delta}}(\mathbf{x})$ is the horizon which is usually taken as a disk or ball of radius $\delta$. $\mathbf{u}(\mathbf{x},t)$ is displacement vector field, and $\mathbf{b}(\mathbf{x},t)$ is a body force density field. $\mathbf{h}(\mathbf{x},t)$ is the prescribed displacement data imposed on the volume constrained boundary $\Omega_c$. 

$\mathbf{C(\x'-\x)}$ is the micromodulus tensor which can be written as\cite{XM}:
\begin{equation}\label{ker:e1}
\mathbf{C}(\mathbf{x^{'}}-\mathbf{x})=\alpha \frac{(\mathbf{x^{'}}-\mathbf{x})\otimes(\mathbf{x^{'}}-\mathbf{x})}{|\mathbf{x^{'}}-\mathbf{x}|^3},
\end{equation}
where $\alpha$ is a scalar parameter introduced to keep the energy of peridynamic model and classic elastic model equal, which is determined by $\delta$ and the elastic modulus $E$:
\begin{equation}
\alpha= \begin{cases}\frac{9E}{\pi\delta^{3}h}& \text {2D plane stress or plane strain},\\\frac{12E}{\pi\delta^{4}}&\text {3D}.\\ \end{cases}
\end{equation}
with $h$ being the plate thickness in the 2D model.

$\mu$ is a history-dependent scalar-valued function which can be written as:
\begin{equation}\label{mu}
\mu(s,t)= \begin{cases} $1$ &\text{if the bond is not broken $s\textless s_0$},\\ $0$ & \text {if the bond is broken $s \ge s_0$}.\end{cases} 
\end{equation}
in which $s$ is the bond stretch, $s_0$ is the critical bond stretch defined by
\begin{equation}\label{fra:s1}
s=\dfrac{|\mathbf{x}^{\prime}+\mathbf{u}^{\prime}-\mathbf{x}+\mathbf{u}|-|\mathbf{x}^{\prime}-\mathbf{x}|}{|\mathbf{x}^{\prime}-\mathbf{x}|},\quad
s_o= \begin{cases}\sqrt{\frac{4 \pi G_o}{9 E \delta}}, & \text { Plane stress }, \\ \sqrt{\frac{5 \pi G_o}{12 E \delta}}, & \text { Plane strain }.\end{cases}
\end{equation}
where $G_0$ is the energy release rate.

The static linear bond-based peridynamic model can be written as:
\begin{equation}\label{mat:e2}
\begin{cases}
-\ds \int_{\mathbf{\mathcal{B}_{\delta}(\mathbf{x})}\cap\Omega} \mu\mathbf{C}(\mathbf{\mathbf{x^{'}}}-\mathbf{x})(\mathbf{u}\left(\mathbf{x}^{\prime}\right)-\mathbf{u}(\mathbf{x})) dV_{\mathbf{\mathbf{x}^{\prime}}}=\mathbf{b}(\mathbf{x})&\text{$\mathbf{x}\in \Omega_s$},\\
\mathbf{u}(\mathbf{x})=\mathbf{h}(\mathbf{x})&\text{$\mathbf{x}\in \Omega_c$}.
\end{cases}
\end{equation}

\subsection{Temporal discretization}
For peridynamic models, several temporal discretization algorithms have bee developed. In this paper, an adaptive dynamic relaxation(ADR) method\cite{KM} is adopted to solve the quasi-static problem Eq.(\ref{mat:e2}) and a second-order Velocity Verlet(VV) algorithm is applied to solve the time-dependent problem Eq.(\ref{mat:e1}).

To solve the quasi-static peridyanmic model, the ADR method is a popular choice which can transform a quasi-static problem into a dynamic problem, which is given by:
\begin{equation}\label{ADR:e1}
\mathbf{D} \ddot{\mathbf{u}}(\mathbf{x}, t)+c\mathbf{D} \dot{\mathbf{u}}(\mathbf{x}, t)=\mathbf{f}+\mathbf{b}
\end{equation}
$\mathbf{D}$ is a fictitious diagonal density matrix. $c$  is a damping coefficient introduced to keep the solution stable. A central-difference explicit integration scheme(CDEI) is used to solve (\ref{ADR:e1}), as is shown in Algorithm \ref{vv1}.
\begin{algorithm}   	
	\SetAlgoNoLine 
	\SetKwInOut{Input}{\textbf{Input}}\SetKwInOut{Output}{\textbf{Output}}   	
	$\dot{\mathbf{u}}^{1/2}=\Delta t \mathbf{D}^{-1} (\mathbf{f}^{0}+\mathbf{b}^{0})/2$ \\	\vspace{0.05in}
	$\dot{\mathbf{u}}^{n+1 / 2}=\left(\left(2-c^{n} \Delta t\right) \dot{\mathbf{u}}^{n-1 / 2}+2 \Delta t \mathbf{D}^{-1} (\mathbf{f}^{n}+\mathbf{b}^{n})\right)/\left(2+c^{n} \Delta t\right)$ \\\vspace{0.05in}
	$\mathbf{u}^{n+1}=\mathbf{u}^{n}+\Delta t \dot{\mathbf{u}}^{n+1 / 2}$
	\caption{\textbf{CDEI}\label{vv1}}
\end{algorithm}

Here,  $\mathbf{f}^{n}$ can be calculated by:
\begin{equation}
\mathbf{\mathbf{f}^{n}}= \int_{\mathbf{\mathcal{B}_{\delta}(\mathbf{x})}\cup\Omega}\mu(s,t^n)\mathbf{C}(\mathbf{\mathbf{x^{'}}}-\mathbf{x})(\mathbf{u}\left(\mathbf{x}^{\prime}, t^{n}\right)-\mathbf{u}(\mathbf{x}, t^{n})) dV_{\mathbf{\mathbf{x}'}}
\end{equation}

For time-dependent problems, VV algorithm is used to integrate Newton's equation of motion, which is described in Algorithm \ref{vv1}. 

\begin{algorithm}
	\SetAlgoNoLine 
	\SetKwInOut{Input}{\textbf{Input}}\SetKwInOut{Output}{\textbf{Output}}    	   	
	$\mathbf{u}^{n+1}=\mathbf{u}^{n}+\Delta t \dot{\mathbf{u}}^n+\Delta t^2(\mathbf{f}^n+\mathbf{b}^n)/2$ \;\vspace{0.05in} 
	$\dot{\mathbf{u}}^{n+1}=\dot{\mathbf{u}}^{n}+\Delta t(\mathbf{f}^n+\mathbf{b}^n+\mathbf{f}^{n+1}+\mathbf{b}^{n+1})/2$\; 
	\caption{\textbf{VV}\label{vv2}}
\end{algorithm}

In conclusion, both the quasi-static and time-dependent problems can be solved by some explicit discretization schemes. The majority of computation cost rests on the calculation of $\mathbf{f}^{n}$.
 \subsection{Spatial discretization}
\begin{figure}[h]
	\centering            
	\includegraphics[scale=0.45]{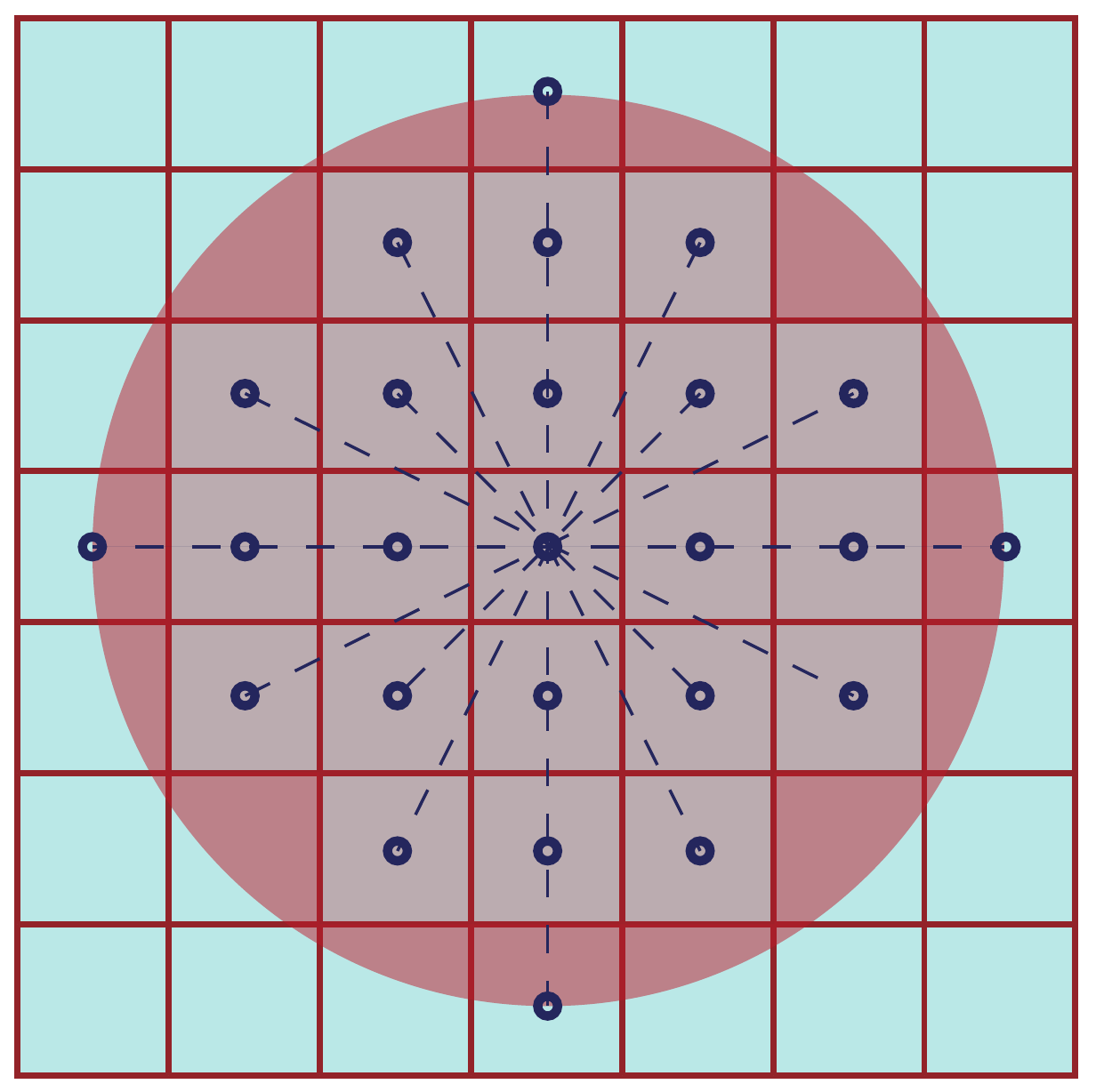}  
	\caption{The volume correction in PD} 
	\label{Fig:1}
\end{figure}

 A meshfree method proposed in \cite{SilAsk} is employed to discretize $\mathbf{f}$ due to its simplicity. The discretized form  of  $\mathbf{f}$ at $\mathbf{x}_p$ can be expressed as:
 \begin{equation}\label{mat:f1}
 \mathbf{f}_{p}=\sum_{\x_q\in \mathbf{\mathcal{B}_{\delta}}(\mathbf{x}_p)\cap \Omega} \mu\mathbf{C}(\mathbf{x}_{q}-\mathbf{x}_{p})(\mathbf{u}_{q}-\mathbf{u}_{p})\lambda_q V_{q},\quad \mathbf{x}\in\Omega_s,
 \end{equation}
where $\mathbf{x}_p$ and $\mathbf{x}_q$ are the positions of material nodes and $V_{q}$ is the volume of node $\mathbf{x}_{q}$. $\lambda_q$ is a volume correction factor which is introduced to correct the volume of neighboring nodes, which are located near the boundary of horizon and partly belong to the horizon, as shown in Fig.\ref {Fig:1}. Lots of volume correction algorithms have been developed. In this paper, $\lambda_q$ is defined as follows\cite{MaOt}:
\begin{equation}\label{fac:e1}
\lambda_{q}=
\begin{cases}
1& \text{when $\|\mathbf{x}_{q}-\mathbf{x}_{p}\|\leq \delta-\dfrac{\Delta x}{2}$} \\
\dfrac{2\delta+\Delta x-2\|\mathbf{x}_{q}-\mathbf{x}_{p}\|}{2\Delta x}&\text{when $\delta-\dfrac{\Delta x}{2}\textless\|\mathbf{x}_{q}-\mathbf{x}_{p}\|\leq \delta$}\\
0&\text{otherwise}
\end{cases}
\end{equation}
where $\Delta x$ is the grid spacing.

For the 2D model, we denote $\mathbf{x}_p=[x_p,y_p]^T$, $\mathbf{x}_{q}=[x_{q},y_{q}]^T$, $\mathbf{u}_p=[u^{x}_{p},u_{p}^{y}]^{T},\mathbf{f}_p=[f_{p}^{x},f_{p}^{y}]^{T}$.
Substituting (\ref{ker:e1}) into (\ref{mat:f1}) yields:
\begin{equation}\label{fxy}
\left[
\begin{matrix}
f_{p}^{x}\\
f_{p}^{x}\\
\end{matrix}
\right]
=
\left[
\begin{matrix}
f_{p}^{xx}+f_{p}^{xy}\\
f_{p}^{yx}+f_{p}^{yy}
\end{matrix}
\right]
\end{equation}
where:
\begin{equation}\label{mod:e2}
\begin{split}
&f_{p}^{xx}=\sum_{\mathbf{x}_{q}\in\mathcal{{B}_{\delta}}(\mathbf{x}_{p})\cap\Omega}\alpha\mu \frac{(x_{q}-x_{p})^{2}}{((x_{q}-x_{p})^{2}+(y_{q}-y_{p})^{2})^{\frac{3}{2}}}\left(u_{q}^{x}-u_{p}^{x}\right)\lambda_q V_{q}\\
&f_{p}^{xy}=\sum_{\mathbf{x}_{q}\in\mathcal{B_{\delta}}(\mathbf{x}_{p})\cap\Omega}\alpha\mu \frac{(x_{q}-x_{p})(y_{q}-y_{p})}{((x_{q}-x_{p})^{2}+(y_{q}-y_{p})^{2})^{\frac{3}{2}}}\left(u_{q}^{y}-u_{p}^{y}\right)\lambda_q V_{q}\\
&f_{p}^{yx}=\sum_{\mathbf{x}_{q}\in\mathcal{B_{\delta}}(\mathbf{x}_{p})\cap\Omega}\alpha\mu \frac{(y_{q}-y_{p})(x_{q}-x_{p})}{((x_{q}-x_{p})^{2}+(y_{q}-y_{p})^{2})^{\frac{3}{2}}}\left(u_{q}^{x}-u_{p}^{x}\right)\lambda_q V_{q}\\
&f_{p}^{xx}=\sum_{\mathbf{x}_{q}\in\mathcal{B_{\delta}}(\mathbf{x}_{p})\cap\Omega}\alpha\mu \frac{(y_{q}-y_{p})^{2}}{((x_{q}-x_{p})^{2}+(y_{q}-y_{p})^{2})^{\frac{3}{2}}}\left(u_{q}^{y}-u_{p}^{y}\right)\lambda_q V_{q}
\end{split}
\end{equation}

 Let $\mathbf{f}^{xx}=[f^{xx}_1,\dots,f^{xx}_N]^T$, $\mathbf{f}^{xy}=[f^{xy}_1,\dots,f^{xy}_N]^T$, $\mathbf{f}^{yx}=[f^{yx}_1,\dots$, $f^{yx}_N]^T,\mathbf{f}^{yy}=[f^{yy}_1,\dots,f^{yy}_N]^T$, $\mathbf{u}^x=[u^x_1,\dots,u^x_N]^T$,  $\mathbf{u}^y=[u^y_1,\dots,u^y_N]^T$ be the $N$-dimension vectors, where $N=N_xN_y$ refer to the number of material points, 
then $\mathbf{f}$ can be rewritten as a matrix-vector multiplication:
\begin{equation}\label{mmv:fm}
\begin{split}
&\mathbf{f}^{xx}=\mathbf{A}^{xx}\mathbf{u}^{x},\mathbf{f}^{xy}=\mathbf{A}^{xy}\mathbf{u}^{y}\\
&\mathbf{f}^{yx}=\mathbf{A}^{xy}\mathbf{u}^{x},\mathbf{f}^{yy}=\mathbf{A}^{yy}\mathbf{u}^{y}
\end{split}   	
\end{equation}

${A}^{xx}$, ${A}^{xy}$, ${A}^{yx}$, ${A}^{yy}$ can be expressed as:
\begin{equation}\label{matrix}
\begin{split}
&A_{p,q}^{xx}=
\begin{cases}
\alpha\mu \dfrac{(x_{q}-x_{p})^{2}}{((x_{q}-x_{p})^{2}+(y_{q}-y_{p})^{2})^{\frac{3}{2}}}\lambda_{q} V_{q}& p\neq q,\mathbf{x}_p \in \Omega_s
\\-\sum\limits_{\mathbf{x}_{q}\in\mathcal{B_{\delta}}(\mathbf{x}_p)\cap\Omega}\alpha\mu \dfrac{(x_{q}-x_{p})^{2}}{((x_{q}-x_{p})^{2}+(y_{q}-y_{p})^{2})^{\frac{3}{2}}}\lambda_{q} V_{q}& p=q,\mathbf{x}_p \in \Omega_s
\end{cases}\\
&A_{p,q}^{xy}=
\begin{cases}
\alpha\mu \dfrac{(x_{q}-x_{p})(y_{q}-y_{p})}{((x_{q}-x_{p})^{2}+(y_{q}-y_{p})^{2})^{\frac{3}{2}}}\lambda_{q} V_{q}& p\neq q,\mathbf{x}_p \in \Omega_s
\\-\sum\limits_{\mathbf{x}_{q}\in\mathcal{B_{\delta}}(\mathbf{x}_p)\cap\Omega}\alpha\mu \dfrac{(x_{q}-x_{p})(y_{q}-y_{p})}{((x_{q}-x_{p})^{2}+(y_{q}-y_{p})^{2})^{\frac{3}{2}}}\lambda_{q} V_{q}& p=q,\mathbf{x}_p \in \Omega_s
\end{cases}\\
&A_{p,q}^{yx}=
\begin{cases}
\alpha\mu \dfrac{(y_{q}-y_{p})(x_{q}-x_{p})}{((x_{q}-x_{p})^{2}+(y_{q}-y_{p})^{2})^{\frac{3}{2}}}\lambda_{q} V_{q}& p\neq q,\mathbf{x}_p \in \Omega_s
\\-\sum\limits_{\mathbf{x}_{q}\in\mathcal{B_{\delta}}(\mathbf{x}_p)\cap\Omega}\alpha\mu \dfrac{(y_{q}-y_{p})(x_{q}-x_{p})}{((x_{q}-x_{p})^{2}+(y_{q}-y_{p})^{2})^{\frac{3}{2}}}\lambda_{q} V_{q}& p=q,\mathbf{x}_p \in \Omega_s
\end{cases}\\
&A_{p,q}^{yy}=
\begin{cases}
\alpha\mu \dfrac{(y_{q}-y_{p})^{2}}{((x_{q}-x_{p})^{2}+(y_{q}-y_{p})^{2})^{\frac{3}{2}}}\lambda_{q} V_{q}& p\neq q,\mathbf{x}_p \in \Omega_s
\\-\sum\limits_{\mathbf{x}_{q}\in\mathcal{B_{\delta}}(\mathbf{x}_p)\cap\Omega}\alpha\mu \dfrac{(y_{q}-y_{p})^{2}}{((x_{q}-x_{p})^{2}+(y_{q}-y_{p})^{2})^{\frac{3}{2}}}\lambda_{q} V_{q}& p=q,\mathbf{x}_p \in \Omega_s
\end{cases}\\
\end{split}
\end{equation}
Since the structure of $\mathbf{A}^{xx}$, $\mathbf{A}^{xy}$, $\mathbf{A}^{yx}$ and $\mathbf{A}^{yy}$ are similar,  we only discuss $\mathbf{f}^{xx}=\mathbf{A}^{xx}\mathbf{u}^x$ in the following sections and record it as $\mathbf{f}=\mathbf{A}\mathbf{u}$ for convenience.
\begin{remark}
As is shown above, the equations in (\ref{matrix}) are satisfied for points $\x_p \in \Omega_s$. The entries in the $p$-th row of matrix $\mathbf{A}$ can be described by (\ref{matrix}). 
If $\x_{p'} \in \Omega_c$, there is no definition for the entries in $p'$-th row of the matrix $\mathbf{A}$.  
In fact, for $\mathbf{x}_p'\in\Omega_c$, the displacement $\mathbf{u}_p'$ is given by the prescribed volume boundary conditions. $A_{p,q}$ in the $p-$th row do not affect the process of getting $\mathbf{u}_p$. Thus $A_{p,q}$ for $\mathbf{x}_p\in\Omega_c$ can be chosen arbitrarily. We will discuss the detailed form of the entries in the $p'$-th row in the following section.
\end{remark}
\section{A fast matrix-based method for the 2D linear bond-based peridynamic model}
\begin{figure}[h]
 	\centering            
 	\includegraphics[scale=0.7]{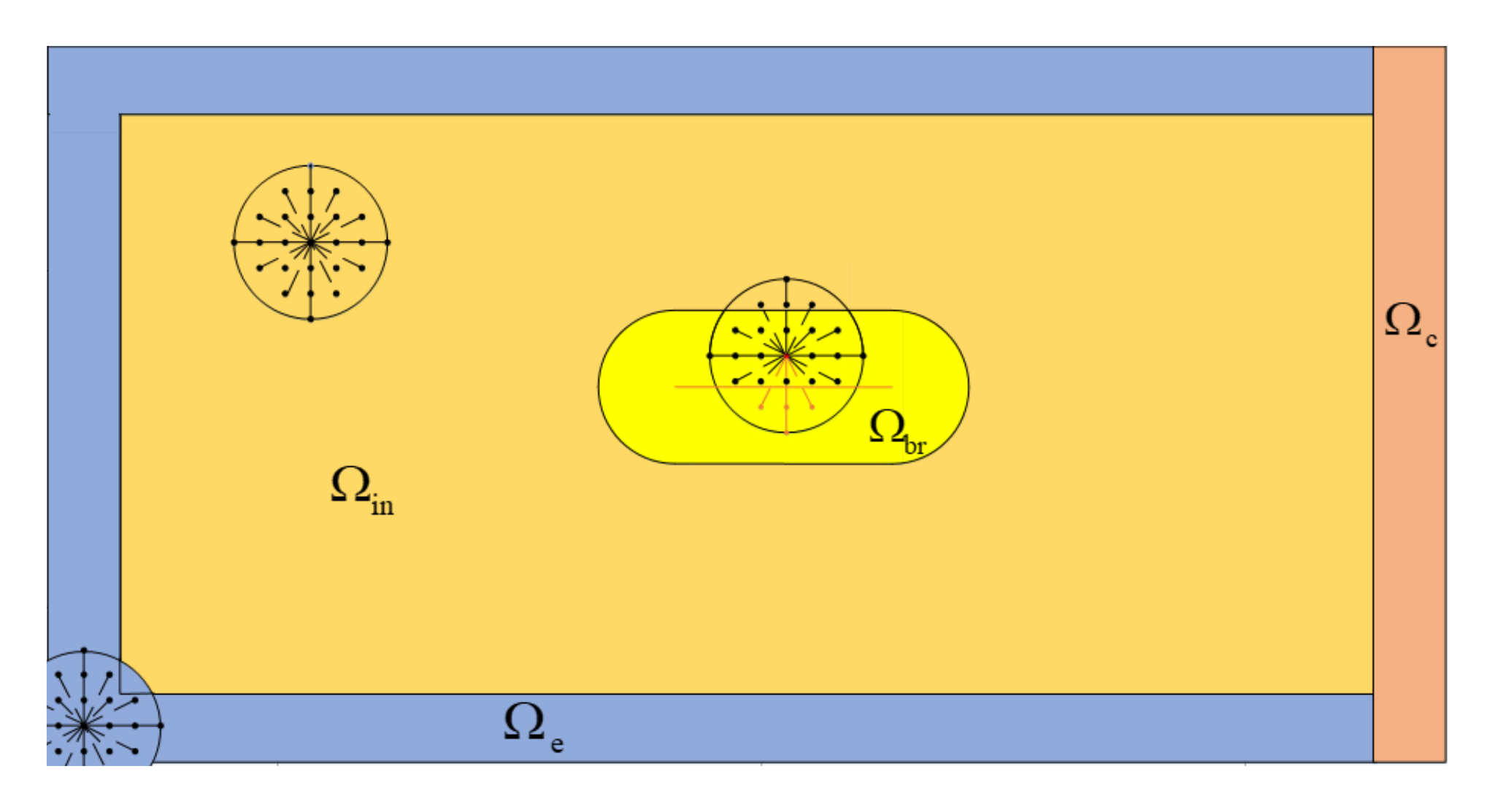}  
 	\caption{Different domains of $\Omega$. The golden part is the sub-domain $\Omega_{in}$. The orange part is the sub-domain $\Omega_c$. The blue part is the sub-domain $\Omega_e$. The yellow part is the sub-domain $\Omega_{br}$.} 
 	\label{Fig:2}
 \end{figure}

 We begin this section by considering a 2D linear bond-based peridynamic model on a rectangular plate $\Omega=\Omega_s\cup\Omega_c=[x_{l},x_{r}]\times[y_{b},y_{t}]$. To investigate the properties of different material points, we divide $\Omega_s$ into $\Omega_{in} \cup\Omega_e \cup\Omega_{br}$, as shown in Fig. \ref{Fig:2}.

 $\bullet$ $\Omega_{in}=\{\mathbf{x}_p: \mathcal{B_{\delta}}(\mathbf{x}_p)\cap\Omega=\mathcal{{B}_{\delta}}(\x_p)$\};
 
 $\bullet$ $\Omega_e=\{\mathbf{x}_p: \mathcal{B_{\delta}}(\mathbf{x}_p)\cap\Omega\neq \mathcal{{B}_{\delta}}(\x_p)$\};
 
  $\bullet$ $\Omega_{br}=\{\mathbf{x}_p: \mathcal{B_{\delta}}(\mathbf{x}_p)\cap\Omega$ \text{include at least one broken bond}\};
  
 Here, $\Omega_{in}$ is the internal area, in which the influence area $\mathcal{B_{\delta}}(\mathbf{x}_p)\cap\Omega$ of all material points on $\Omega_{in}$ is a complete disk $\mathcal{B}_{\delta}(\mathbf{x}_p)$. The sub-domain $\Omega_e$ includes the points which are still in $\Omega_s$, but so close to the boundary of $\Omega$ that the influence area $\mathcal{B_{\delta}}(\mathbf{x}_p)\cap\Omega$ is no longer a complete disk $\mathcal{B}_{\delta}(\mathbf{x}_p)$. There is at least one material point $\mathbf{x}_q$ in the influence area of material point $\mathbf{x}_p\in\Omega_{br}$, which causes the bond between $\mathbf{x}_p$ and $\mathbf{x}_q$ to be broken. If $\Omega_e=\Omega_{br}=\varnothing$, in other words, the domain $\Omega$ only contains $\Omega_{in}$ and $\Omega_{c}$, this model will turn into the case mentioned in \cite{XH}.

 The domain is discretized with a uniform spatial partition. $N_{x}$, $N_{y}$ denote the numbers of nodes in the $x$-direction and $y$-direction. $\mathbf{x}_{p}$ can be expressed as $\mathbf{x}_{p}=[x_{i},y_{j}]^{T}$,  where $x_{i}=x_{l}+(i-1/2)h_{x}$  and $y_{j}=y_{b}+(j-1/2)h_{y}$. $h_x$ and $h_y$ are the grid spacing in each direction. Without loss of generality, we let $h_x=h_y=h$.

 The node number $p$ mentioned in (\ref{mat:f1}) is 
 related to $i$ and $j$ by $p=(j-1)N_x+i$.  In order to discuss the matrix structure more conveniently, $f_p$ and $u_p$ in (\ref{mod:e2}) are rewritten as a $f_{i,j},u_{i,j}$ respectively.Then $\mathbf{f}$ and $\mathbf{u}$ can be rewritten as:
\begin{equation}
\begin{split}
&\mathbf{u}=[u_{1,1},\dots,u_{N_x,1},\dots,u_{1,N_{y}},\dots,u_{N_{x},N_{y}}]^{T},\\
&\mathbf{f}=[f_{1,1},\dots,f_{N_x,1},\dots,f_{1,N_y},\dots,f_{N_{x},N_{y}}]^{T},
\end{split}
\end{equation}
in which $f_{i,j}$ can be reorganized as follows for $\mathbf{x}_q=(x_{i^{'}},y_{j^{'}})$:
\begin{equation}\label{fce:e1}
f_{i,j}=\sum_{\mathbf{x}_{q}\in\mathcal{B_{\delta}}(\x_p)\cap\Omega}\alpha\mu \frac{(x_{i^{'}}-x_{i})^{2}}{((x_{q}-x_{p})^{2}+(y_{q}-y_{p})^{2})^{\frac{3}{2}}}\left(u_{i^{'},j^{'}}-u_{i,j}\right)\lambda_{i^{'},j^{'}}V_{i^{'},j^{'}}\\
\end{equation}
With the uniform node distribution, the volume $V_{i^{'},j^{'}}$ of each material point is $h^2$.

\subsection{The block-banded-block structure of matrix $\mathbf{A}$}
For general regions, the matrix has no specific structure. But under the uniform mesh of the rectangular domain, matrix $\mathbf{A}$ can be rewritten into such a form as follows according to Eq.(\ref{fce:e1}):

\begin{equation}\label{stu:bad}
\mathbf{A}=
\left[
\begin{matrix}
\mathbf{B_{1,1}} &\cdots & \mathbf{B_{1,N_{y}}} \\
\vdots & \ddots &\vdots\\
\mathbf{B_{N_{y},1}} & \cdots & \mathbf{B_{N_{y},N_{y}}}
\end{matrix}
\right],   \mathbf{B}_{i,i^{'}}=
\left[
\begin{matrix}
c^{1,1}_{i,i^{'}} &  \cdots & c^{1,N_{x}}_{i,i^{'}} \\
\vdots &\ddots & \vdots \\
c^{N_{x},1}_{i,i^{'}} & \cdots & c^{N_{x},N_{x}}_{i,i^{'}}
\end{matrix}
\right],
\end{equation}
where $A_{p,q}=c_{i,i^{'}}^{j,j^{'}}$ can be expressed as:
\begin{equation}\label{cij}
	c_{i,i^{'}}^{j,j^{'}}=
	\begin{cases}
	\alpha\mu \dfrac{(y_{j^{'}}-y_{j})(x_{i^{'}}-x_{i})}{((x_{q}-x_{p})^{2}+(y_{q}-y_{p})^{2})^{\frac{3}{2}}}\lambda_{i^{'},j^{'}}h^2& q\neq p,\mathbf{x}_p \in \Omega_s
	\\-\sum\limits_{\mathbf{x}_{q}\in\mathcal{B_{\delta}}(\mathbf{x}_p)\cap\Omega}\alpha\mu \dfrac{(y_{j^{'}}-y_{j})(x_{i^{'}}-x_{i})}{((x_{q}-x_{p})^{2}+(y_{q}-y_{p})^{2})^{\frac{3}{2}}}\lambda_{i^{'},j^{'}}h^2& q=p,\mathbf{x}_p \in \Omega_s
	\end{cases}
\end{equation}
$A_{p,q}$ represents the action of the material point $\mathbf{x}_q=[x_{i^{'}},y_{j^{'}}]^T$ on $\mathbf{x}_p=[x_i,y_j]^T$, and the matrix block $\mathbf{B}_{i,i^{'}}$ represents the action of  material points in $i^{'}-$th row on material points in $i-$th row. Furthermore, each entry $c^{j,j^{'}}_{i,i^{'}}$ in the matrix block $\mathbf{B}_{i,i^{'}}$ represents the effect of the $j^{'}$-th material point on the $i^{'}$-th row  on the $j$-th point on the $i$-th row.

In the PD model, only the material point $\mathbf{x}_q\in\mathcal{B}_{\delta}(\mathbf{x}_p)\cap\Omega$ interacts with the material point $\mathbf{x}_p$, which means:
\begin{equation}\label{band}
A_{p,q}=c^{j,j^{'}}_{i,i^{'}}=0,\text{ if $(x_{i^{\prime}}-x_{i})^{2}+(y_{j^{\prime}}-y_{j})^{2}> \delta^2$}
\end{equation}
It's easy to see that $c^{j,j^{'}}_{i,i^{'}}=0$ if $|j^{'}-j|>M$ or $|i^{'}-i|>M$, $M=\delta/h$. Each matrix block satisfies $\mathbf{B}_{i,i^{'}}=0$ if $|i^{'}-i|>M$. Furthermore, if $|i^{'}-i|\leq M$ and $|j^{'}-j|>M$, the entries in the matrix block $\mathbf{B}_{i,i^{'}}$ also satisfies $c^{j,j^{'}}_{i,i^{'}}=0$.  Therefore $\mathbf{A}$ is a block-banded-block matrix, as shown in Fig. {\ref{Fig:3}.
		\begin{figure}[h]
		\centering            
		\includegraphics[scale=0.5]{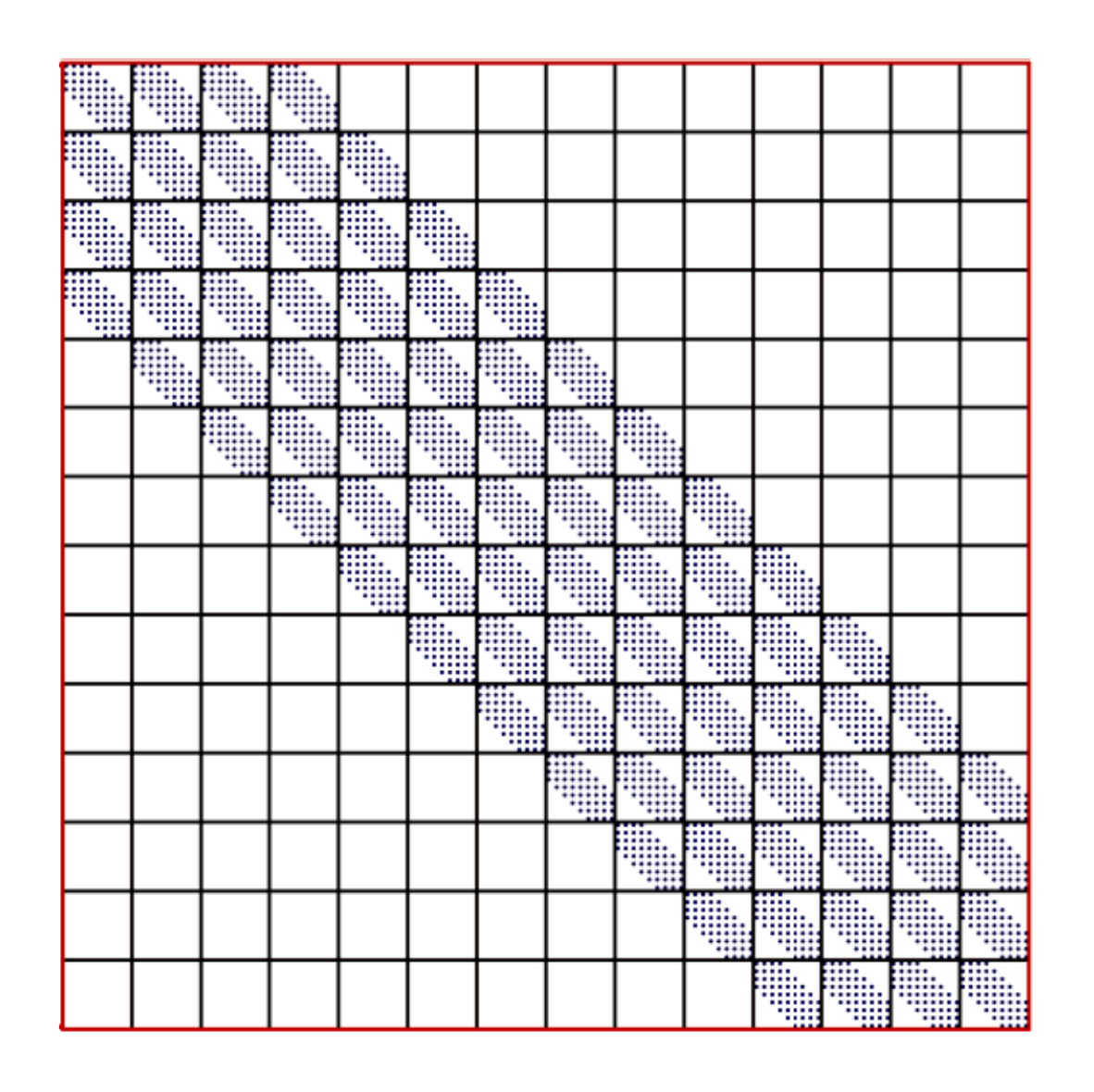}  
		\caption{The block-banded-block structure of matrix $\mathbf{A}$ for 2D PD} 
		\label{Fig:3}
	\end{figure} 
\subsection{Analysis of $A_{p,q}$ for any $\mathbf{x}_p\in\Omega_{in}$} 
In most cases, $\Omega_{in}$ occupies most part of the $\Omega$. Different from the material points in $\Omega_e$ and $\Omega_{br}$, The influence area  of these material points  in $\Omega_{in}$ is a complete disk , in which none of the bonds  inside $\Omega_{in}$ are broken.

If $\mathbf{x}_p\in\Omega_{in}$, we have $\mathcal{B_{\delta}}(\mathbf{x}_p)\cap\Omega=\mathcal{B_{\delta}}(\mathbf{x}_p)$, $\mu=1$, then $f_{i,j}$ can be expressed as follows:
\begin{equation}\label{fce:e2}
f_{i,j}=\sum_{\mathbf{x}_{q}\in\mathcal{B_{\delta}}(\mathbf{x}_p)}\alpha \frac{(x_{i^{'}}-x_{i})^{2}}{((x_{q}-x_{p})^{2}+(y_{q}-y_{p})^{2})^{\frac{3}{2}}}\left(u_{i^{'},j^{'}}-u_{i,j}\right)\lambda_{i^{'},j^{'}}h^2, \quad \forall \x_p\in\Omega_{in}\\
\end{equation} 
    
 With the relations $x_{i}=x_{l}+(i-1/2)h$ and $y_{j}=y_{b}+(j-1/2)h$, matrix entries $A_{p,q}$ can be reorganized as follows when $\mathbf{x}_p\in\Omega_{in}$ and $q\neq p$:
  \begin{equation}\label{frb:e1}
  \begin{split} 
  	A_{p,q}=c_{i,i^{'}}^{j,j^{'}}&=\alpha\frac{(x_{i^{\prime}}-x_{i})^{2}}{((x_{i^{'}}-x_{i})^{2}+(y_{i^{'}}-y_{i})^{2})^{\frac{3}{2}}}\lambda_{i^{'},j^{'}}h^2\\
  	&=\alpha\frac{(i^{\prime}-i)^{2}}{((i^{\prime}-i)^{2}+(j^{\prime}-j)^{2})^{\frac{3}{2}}}\lambda_{i^{'},j^{'}}h.
  \end{split}
  \end{equation}
    where
  \begin{equation}\label{tra:e2}  
  \lambda_{i^{'},j^{'}}=
  \begin{cases}
  1&\text{when $\|\mathbf{x}_{{q}}-\mathbf{x}_{p}\|\leq \delta-h/2$} \\
  M-\sqrt{(i^{\prime}-i)^{2}+(j^{\prime}-j)^{2}}/2+1/2&\text{when $\delta-h/2\textless\|\mathbf{x}_{q}-\mathbf{x}_{p}\|\leq \delta$}\\
  0&\text{otherwise}
  \end{cases}
  \end{equation}
  When $q=p$, we can also get:
  \begin{equation}\label{frb:e4}
  A_{p,p}=c_{i,i}^{j,j}=-\sum\limits_{\mathbf{x}_{q}\in\mathcal{B_{\delta}}(\mathbf{x}_p)}c_{i,i^{'}}^{j,j^{'}}=-\sum\limits_{\mathbf{x}_{q}\in\mathcal{B_{\delta}}(\mathbf{x}_p)}\alpha\frac{(i^{\prime}-i)^{2}}{((i^{\prime}-i)^{2}+(j^{\prime}-j)^{2})^{\frac{3}{2}}}\lambda_{i^{'},j^{'}}h
  \end{equation}
    Thus the entries $A_{p,q}$ does not depend on the position of $\mathbf{x}_p$ or $\mathbf{x}_q$, but on the distance $\mathbf{x}_p-\mathbf{x}_q$. Moreover with a uniform mesh on a rectangular plate, $A_{p,q}$ is directly related to the difference $i^{'}-i$ and $j^{'}-j$, according to Eq.(\ref{frb:e4}). Since $-M\leq i^{'}-i,j^{'}-j \leq M$, for every $p-$ th row if $\mathbf{x}_p\in\Omega_{in}$, there are $(2M+1)^2$ non-zero entries and they are the same . Therefore, we only need to store these matrix entries in a $2M+1$-by-$2M+1$ matrix  instead of a $N$-by-$(2M+1)^2$ matrix, which greatly reduce memory size. and footprint. The $2M+1$-by-$2M+1$ matrix is termed as the kernel matrix, which can be defined as:
    \begin{equation}\label{frb:e2}
    K_{m,n}=K_{i^{'}-i+M+1,j^{'}-j+M+1}=c_{i,i^{'}}^{j,j^{'}},\quad 0\leq m,n \leq 2M+1,\quad \forall \x_p\in\Omega_{in}
    \end{equation}
    In the actual calculation, we can obtain $K_{m,n}$ by computing the interaction between $\mathbf{x}_p$ and $(2M+1)^2$ material points around. Here $\mathbf{x}_p$ is an arbitrary material point in $\Omega_{in}$, as shown in Fig. \ref{Fig:22}. Specifically, $K_{M+1,M+1}=c_{i,i}^{j,j}=A_{p,p}$.

    \begin{figure}[h]
	\centering    	
	\subfigure[]
	{
		\begin{minipage}{6cm}
			\centering          
			\includegraphics[scale=0.5]{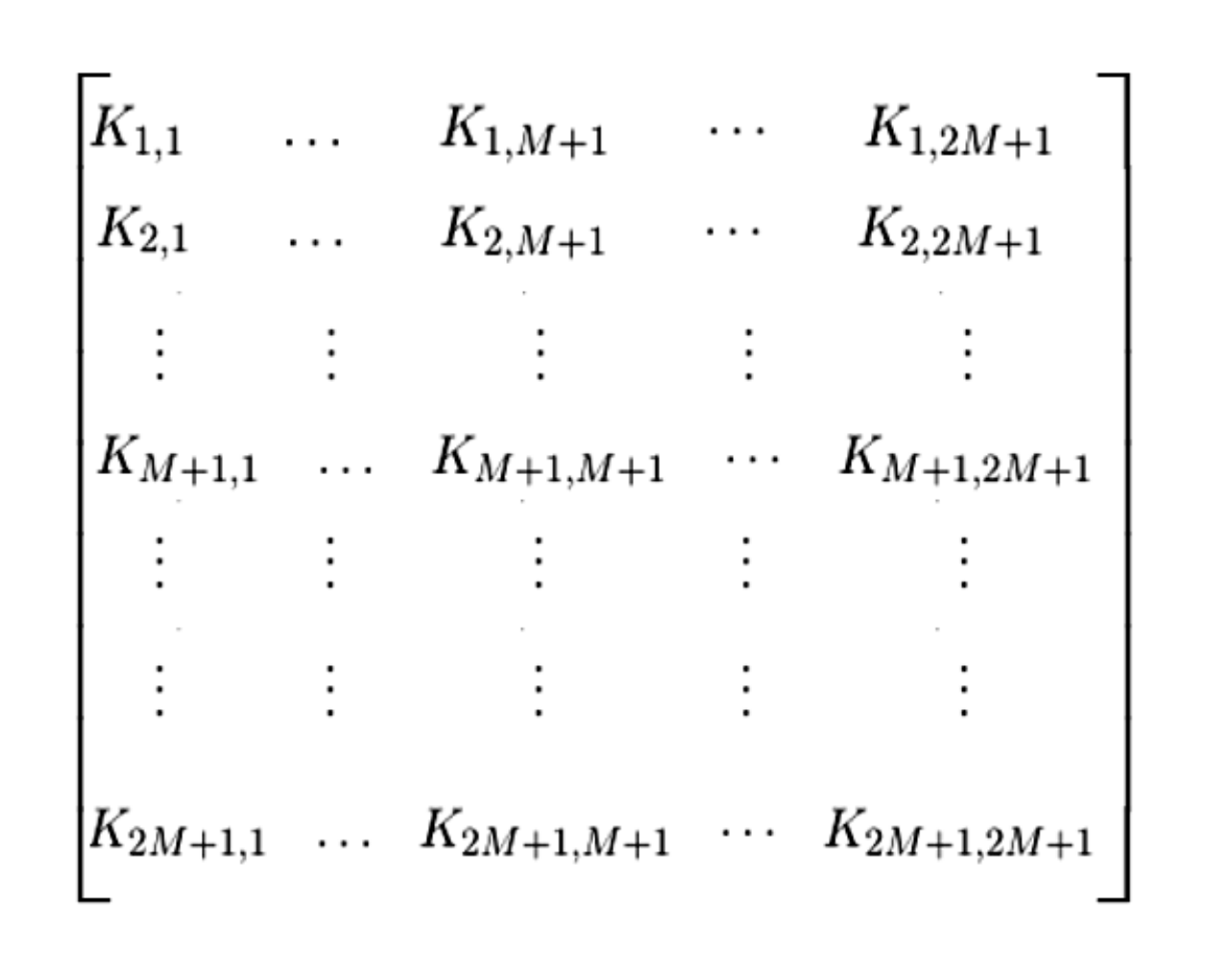}  
		\end{minipage}
	}	
	\subfigure[]
	{
		\begin{minipage}{6cm}
			\centering      
			\includegraphics[scale=0.5]{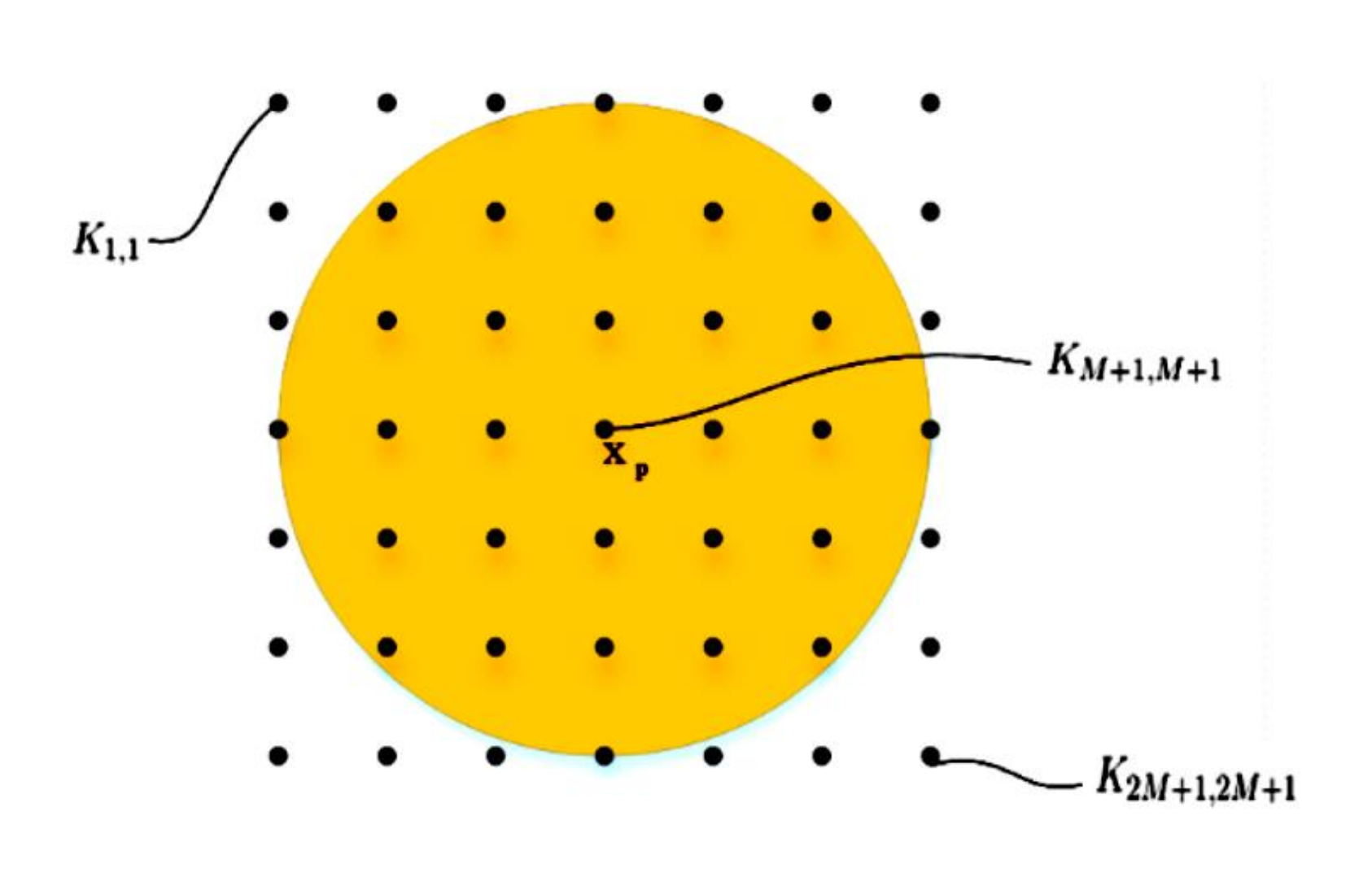}   
		\end{minipage}
	}	
	\caption{Description of the keneral matrix: (a) Entries $K_{m,n}$ in the kernel matrix (b) Contact between material points and matrix entries}
	\label{Fig:22}
\end{figure}

 Assuming that all material points are in $\Omega_{in}$, which means all material points satisfy the Eqs.(\ref{frb:e1}) and (\ref{frb:e4}), we can get a matrix $\hat{\mathbf{A}}$ defined by:
 \begin{equation}
 	\mathbf{\hat{A}}=
 	\left[
 	\begin{matrix}
 	\mathbf{\hat{B}_{1,1}} &\cdots & \mathbf{\hat{B}_{1,N_{y}}} \\
 	\vdots & \ddots &\vdots\\
 	\mathbf{\hat{B}_{N_{y},1}} & \cdots & \mathbf{\hat{B}_{N_{y},N_{y}}}
 	\end{matrix}
 	\right],   \mathbf{\hat{B}}_{i,i^{'}}=
 	\left[
 	\begin{matrix}
 	\hat{c}^{1,1}_{i,i^{'}} &  \cdots & \hat{c}^{1,N_{x}}_{i,i^{'}} \\
 	\vdots &\ddots & \vdots \\
 	\hat{c}^{N_{x},1}_{i,i^{'}} & \cdots & \hat{c}^{N_{x},N_{x}}_{i,i^{'}}
 	\end{matrix}
 	\right],
 \end{equation}
 This is a block-banded-block matrix generated by $\mathbf{K}$, which means:
   \begin{equation}\label{frb:e2}
    \hat{A}_{p,q}=\hat{c}_{i,i^{'}}^{j,j^{'}}=K_{i^{'}-i+M+1,j^{'}-j+M+1}=K_{m,n},\quad \forall \mathbf{x}_p\in\Omega
   \end{equation}
   
The value of $\hat{A}_{p,q}$ depends only on $i^{'}-i$ and $j^{'}-j$. Hence, for every matrix block $\mathbf{\hat{B}}_{i,i^{'}}$, entries $\hat{c}_{i,i^{'}}^{j,j^{'}}$ on each diagonal  are equal. For the matrix $\mathbf{\hat{A}}$,  blocks $\mathbf{\hat{B}}_{i,i^{'}}$ on each diagonal are equal. A matrix satisfies the above properties  is called a Toeplitz-Block-Toeplitz(TBT) matrix.

To construct the fast method , we decompose the matrix $\mathbf{A}=\mathbf{\hat{A}+(A-\hat{A})}$. Thus we have:  
  \begin{equation}\label{dec:e1}
  \mathbf{f}=\mathbf{f^{'}}+(\mathbf{A}-\hat{\mathbf{A}})\mathbf{u},
  \end{equation} 
  where $\mathbf{f^{'}}=\mathbf{\hat{A}}\mathbf{u}$.  Based on the TBT structure of $\hat{\mathbf{A}}$, The matrix-vector multiplication $\mathbf{\hat{A}}\mathbf{u}$ can be accelerated by a fast matrix-vector multiplication(FMVM)\cite{XH} in Algorithm \ref{vv3}:

  \begin{algorithm}   	
  	\SetAlgoNoLine 
  	\SetKwInOut{Input}{\textbf{Input}}\SetKwInOut{Output}{\textbf{Output}}   	
  	$\mathbf{\hat{G}}=\mathbf{FFT2}(\mathbf{G}),\mathbf{\hat{U}}=\mathbf{FFT2}(\mathbf{U}$) \\	\vspace{0.03in}
  	$\mathbf{\hat{H}}=\mathbf{\hat{G}}\circ\mathbf{\hat{U}}$ \\ \vspace{0.03in}
  	$\mathbf{H}=\mathbf{FFT2^{-1}}(\mathbf{\hat{H}})$\\\vspace{0.03in}
  	\caption{\textbf{FMVM}\label{vv3}}
  \end{algorithm}

  Here we use $\mathbf{FFT2}$ and $\mathbf{FFT2^{-1}}$ to denote two-dimension FFT and iFFT operations. $\mathbf{\hat{H}}$ represents the Hadamard product of $\mathbf{\hat{G}}$ and $\mathbf{\hat{U}}$. $\mathbf{H}$ is a $2N_x$-by-$2N_y$ matrix and we can obtain $\mathbf{f}$ by $f_{jN_x+i}=H_{i,j}$, if $i\leq N_x$ and  $j\leq N_y$. 
  $\mathbf{G}$ is the first column of a extended matrix embedded by $\mathbf{K}$. Here we write it as a $2N_x$-by-$2N_y$ matrix , namely:
  \begin{equation}
  \mathbf{G}:=
  \left(\begin{smallmatrix}
  K_{M+1,M+1} & \cdots & K_{M+1, 2M+1} & 0 & \cdots & 0 & K_{M+1,1} & \cdots & K_{M+1,M} \\
  \vdots & \ddots & \vdots & \vdots & \ddots & \vdots & \vdots & \ddots & \vdots \\
  K_{2M+1, M+1} & \cdots & K_{2M+1, 2M+1} & 0 & \ddots & 0 & K_{2M+1,1} & \cdots &K_{2M+1,M} \\
  0 & \cdots &  0& 0& \cdots & 0 & 0& \cdots & 0 \\
  \vdots & \ddots &\ddots & \ddots & \ddots & \ddots & \ddots & \ddots & \vdots \\
  0 & \cdots & 0 & 0 & \cdots & 0 & 0 & \cdots & 0 \\
  K_{1, M+1} & \cdots & K_{1, 2M+1} & 0 & \ddots & 0 & K_{1,1} & \cdots & K_{2M+1,M} \\
  \vdots & \ddots & \vdots & \vdots & \ddots & \vdots & \vdots & \ddots & \vdots \\
  K_{M,M+1} & \cdots & K_{M, 2M+1}& 0 & \cdots & 0 & K_{M,1} & \cdots & K_{M,M}
  \end{smallmatrix}\right)
  \end{equation}
  which means each matrix entries $G_{i,j}$ can be expressed as follows:  
  \begin{equation}
  G_{i,j}=
  \begin{cases}
  0&\text{ $i\in[M+2,2N_x-M] \text{ or j}\in[M+2,2N_y-M]$}\\
  K_{m,n}&\text{otherwise}\\
  \end{cases}  	
  \end{equation}
  where 
  \begin{equation}
  \begin{split}
  &m=\begin{cases}
  i+M&  i\in[1,M+1]\\
  i-2N_x+M&  i\in[2N_x-M+1,2N_x] 
  \end{cases}\\
  &n=\begin{cases}
  j+M&  j\in[1,M+1]\\
  j-2N_y+M& j\in[2N_y-M+1,2N_y]
  \end{cases}
  \end{split}
  \end{equation} 
 $\mathbf{U}$ is  a extended matrix embedded by displacement vector $\mathbf{u}$, which can be expressed as:
  \begin{equation}
  \mathbf{U}:=
  \left(\begin{matrix}
  u_{1,1} & \cdots & u_{N_x,1} & 0 & \cdots & 0 \\
  u_{1,2} & \cdots & u_{N_x,2} & 0 & \cdots & 0 \\
  \vdots&\vdots&\vdots&\vdots&\vdots&\vdots\\
  u_{1,N_y-1} & \cdots & u_{N_x,N_y-1} & 0 & \cdots & 0 \\
  0 & \cdots & 0 & 0& \cdots & 0 \\
  \vdots&\vdots&\vdots&\vdots&\vdots&\vdots\\
   0 & \cdots & 0 & 0& \cdots & 0 \\

  \end{matrix}\right)
  \end{equation}
  which means each entries $U_{i,j}$ can be defined by:
  \begin{equation}
  	U_{i,j}=
  	\begin{cases}
  	u_{i,j}&\text{if $i\leq N_x$ and $j\leq N_y$} \\
  	0&\text{otherwise}
  	\end{cases}
  \end{equation}
  
  By Algorithm \ref{vv3}, The calculation  of form $\mathbf{f^{'}}=\mathbf{\hat{A}u}$ can be decreased from $O(N^2)$ to $O(N\log N)$. The calculation of $\mathbf{(A-\hat{A})u}$ needs to be considered additionally, which will be discussed below.

  	Algorithm \ref{vv3} is implemented with the Matlab code. The codes of the two-dimensional Fourier transform and the two-dimensional inverse transformation are called as $\hat{\mathbf{G}}=fft2(\mathbf{G})$ and ${\mathbf{H}}=ifft2(\mathbf{\hat{H}})$.

  \subsection{Analysis of $A_{p,q}$ for any $\mathbf{x}_p\in\Omega_c$}
In many applications it is usually desired/needed to apply local boundary conditions, the properties of material points $\mathbf{x}_p\in\Omega_c$ are not considered in the matrix, so that the actual matrix is not a square matrix. However, the FMVM algorithm requires that the matrix be a square matrix, so we need to consider it in the form $\mathbf{f=Au}$.

  According to Eq.(\ref{mat:e1}) and Eq.(\ref{mat:e2}), $u_p$ can be expressed as:
  \begin{equation}\label{up}
  	u_p=h(\mathbf{x}_p), \quad \mathbf{x}_p\in\Omega_c
 \end{equation}  
	Equation (\ref{up}) means $u_p$ is given by the prescribed displacement data rather than the form $\mathbf{f=Au}$. Thus $A_{p,q}$ for $\mathbf{x}_p\in\Omega_c$ is meaningless and can be arbitrarily chosen. Here we let $A_{p,q}=\hat{A}_{p,q}$, which means $A_{p,q}-\hat{A}_{p,q}=0$, so that these entries do not have to repeat operations in the form $(\mathbf{A}-\hat{\mathbf{A}})\mathbf{u}$.
. When the computation domain $\Omega=\Omega_c\cup\Omega_{in}$,  the stiff matrix $\mathbf{A}=\mathbf{\hat{A}}$. 

Considering that the displacement $u_p$ is obtained through $h(\mathbf{x}_p)$ for $\mathbf{x}_p\in\Omega_c$ , we need to replace the displacement with $h(\mathbf{x}_p)$ after obtaining the displacement through $\mathbf{f}$ in each time step of Algorithm \ref{vv1} or Algorithm \ref{vv2}. This will bring additional calculation, but it will not exceed $O(N\log N)$ in general.

\subsection{Analysis of $A_{p,q}$ for any $\mathbf{x}_p\in\Omega_e$}
For the material points $\mathbf{x}_p\in\Omega_e$, the influence area  is not a complete disk. Hence, according to (\ref{cij}), if $\mathbf{x}_p\in\Omega_e$, most of the matrix entries in $p-$th row are equal for the matrix entries corresponding to the material points on $\Omega_{in}$ except the matrix entries on the main diagonal.

When $\mathbf{x}_p\in\Omega_e$, we have $\mathcal{B_{\delta}}(\mathbf{x}_p)\cap\Omega\neq\mathcal{B_{\delta}}(\mathbf{x}_p)$ and $\mu=1$. Thus $f_{p}$ can be written as:
\begin{equation}
f_{i,j}=\sum_{\mathbf{x}_{q}\in\mathcal{B_{\delta}}(\mathbf{x}_p)\cap\Omega}\alpha \frac{(x_{i^{'}}-x_{i})^{2}}{((x_{i^{'}}-x_{i})^{2}+(y_{i^{'}}-y_{i})^{2})^{\frac{3}{2}}}\lambda_{i^{'},j^{'}}h^2\left(u_{i^{'},j^{'}}-u_{i,j}\right)\lambda_{i^{'},j^{'}}h^2, \quad \forall \x_p\in\Omega_{e}\\
\end{equation} 
Then entries $A_{p,q}$ can be represented as:
\begin{equation}
A_{p,q}=
\begin{cases}
\alpha \dfrac{(y_{j^{'}}-y_{j})(x_{i^{'}}-x_{i})}{((x_{i^{'}}-x_{i})^{2}+(y_{i^{'}}-y_{i})^{2})^{\frac{3}{2}}}\lambda_{i^{'},j^{'}}h^2& q\neq p
\\-\sum\limits_{\mathbf{x}_{q}\in\mathcal{B_{\delta}}(\mathbf{x}_p)\cap\Omega}\alpha \dfrac{(y_{j^{'}}-y_{j})(x_{i^{'}}-x_{i})}{((x_{i^{'}}-x_{i})^{2}+(y_{i^{'}}-y_{i})^{2})^{\frac{3}{2}}}\lambda_{i^{'},j^{'}}h^2& q=p
\end{cases}
\end{equation}
Notice that each entry in $p-$th is equal to that in $q-$th row expect one on the diagonal for $\mathbf{x}_p\in\Omega_e$ and $\mathbf{x}_q\in\Omega_{in}$, which means:
\begin{equation}	
\begin{cases}
   &A_{p,q}=\hat{A}_{p,q},\quad \quad q\neq p\\
   &A_{p,p}=-\sum\limits_{\mathbf{x}_{q}\in\mathcal{B}_{\delta}(\mathbf{x}_p)\cap\Omega} A_{p,q}=-\sum\limits_{\mathbf{x}_{q}\in\mathcal{B}_{\delta}(\mathbf{x}_p)\cap\Omega} \hat{A}_{p,q}
	\neq-\sum\limits_{\mathbf{x}_{q}\in\mathcal{B}_{\delta}(\mathbf{x}_p)} \hat{A}_{p,q}=\hat{A}_{p,p}\\
\end{cases}
\end{equation}

Hence, $A_{p,p}$ can be expressed as follows:
\begin{equation}\label{Oe:e2}
	A_{p,p}=\hat{A}_{p,p}+\sum\limits_{\mathbf{x}_q\in\mathcal{B}_{\delta}(\mathbf{x}_p\backslash(\mathcal{B}_{\delta}(\mathbf{x}_p)\cap\Omega)}A_{p,q}
\end{equation}
 Then the matrix $\mathbf{A}$ can be decomposed into following form by introducing a diagonal matrix $\mathbf{D}$:
\begin{equation}\label{seb:e2}
\mathbf{A}=\mathbf{\hat{A}}+\mathbf{D}^e
\end{equation}
Here matrix entries in $\mathbf{D}^e$ can be written as:
\begin{equation}\label{mat:se}
D_{p,p}^e=
\begin{cases}
\sum\limits_{\mathbf{x}_q\in\mathcal{B}_{\delta}(\mathbf{x}_p)\backslash(\mathcal{B}_{\delta}(\mathbf{x}_p)\cap\Omega)}A_{p,q}&\text{if $\mathbf{x}_p\in \Omega_{e}$}
\\0& \text{otherwise}
\end{cases}
\end{equation}

Then $\mathbf{f}$ can be decomposed as follows:
\begin{equation}\label{tea:sc}
f_{p}=
\begin{cases}
f^{'}_p+D_{p,p}^{e}u_{p}&\text{if $\mathbf{x}_p\in \Omega_{e}$}
\\f^{'}_p & \text{otherwise}
\end{cases}
\end{equation}
The total number of material points on $\Omega_e$ do not exceed the total number of material points $N$. Thus the calculation and storage memory brought by (\ref{tea:sc}) is $O(N)$.

Since the displacement of the material point on $\Omega_c$ is not affected by the form $\mathbf{Au=f}$, the problem of incomplete horizon on $\Omega_c$ do not need to be considered, which means that $\Omega_e\cap \Omega_c=\varnothing$ in most cases. 

\begin{remark}
	A special case is that a material point is affected by two constraints, which means: (a) The displacement constraint conditions is applied in the $x$-direction, so it is regarded as a material point on $\Omega_c$, and the displacement $u_x$ is replaced when calculating $\mathbf{f^x=A^{xx}u^x+A^{xy}u^y}$; (b) It is affected by the incomplete horizon in the $y$-direction, so it is considered as a matter point on $\Omega_e$ and the corresponding $\mathbf{Du}$ is subtracted when calculating $\mathbf{A^{yx}u^{y}}, \mathbf{A^{yy}u^{y}}$.
\end{remark}
\begin{remark}
Surface correction algorithms are often used on material points $\mathbf{x}_p\in\Omega_e$ to calculate accurately in engineering problems\cite{QF}. In this algorithm, a coefficient $v_{p,q}$ is introduced to increase the micromodule of each bond in $\mathcal{B}_{\delta}(\mathbf{x}_p)\cap\Omega$, which means $A_{p,q}=v_{p,q}\hat{A}_{p,q}$ for $\mathbf{x}_p\in\Omega_e$ and $p\neq q$. This way, the impact that $\mathcal{B}_{\delta}(\mathbf{x}_p)\cap\Omega$ is not a complete disk is eliminated. However, this coefficient breaks the above matrix structure, so the form $A_{p,q}u_{q}$ for $\mathbf{x}_p\in\Omega_e$ needs to be recalculated. We will mention this part in numerical examples. 
\end{remark}

\subsection{Analysis of $A_{p,q}$ for $\mathbf{x}_p\in\Omega_{br}$}
 For the sub-domain $\Omega_{br}$, the interaction between the material points on the broken bond is considered as $0$. Hence,  the structure of the matrix mentioned above is destroyed, which requires special treatment.

 Since all sub-domains with cracks are called $\Omega_{br}$ for the domain $\Omega$, $\Omega_{br}$ does not exist alone. The sub-domain always intersects one of  sub-domains $\Omega_{in}$, $\Omega_{e}$, and $\Omega_{c}$, which means $\Omega_{br}\cap(\Omega_s\cup\Omega_c)\neq \varnothing$. In this case, we first consider the effect of broken bonds on the matrix. For $\mathbf{x}_p\in\Omega_{br}$, $f_p$ can be expressed as follows:
\begin{equation}
f_{i,j}=\sum_{\mathbf{x}_{q}\in\mathcal{B}_{\delta}(\mathbf{x}_p)\cap\Omega}\alpha\mu \frac{(x_{i^{'}}-x_{i})^{2}}{((x_{i^{'}}-x_{i})^{2}+(y_{i^{'}}-y_{i})^{2})^{\frac{3}{2}}}\left(u_{i^{'},j^{'}}-u_{i,j}\right)\lambda_{i^{'},j^{'}}h^2\\
\end{equation}
Here the history-dependent scalar valued function $\mu$ is defined in (\ref{mu}).

As observed from Fig. \ref{Fig:4}, when the bond between two material points $\x_p$ and $\x_q$ is broken, the history-dependent scalar valued function $\mu=0$. Thus $A_{p,q}=0\neq\hat{A}_{p,q}$ for $\mathbf{x}_p\in\Omega_{br}$.
\begin{figure}[h]
	\centering            
	\includegraphics[scale=0.45]{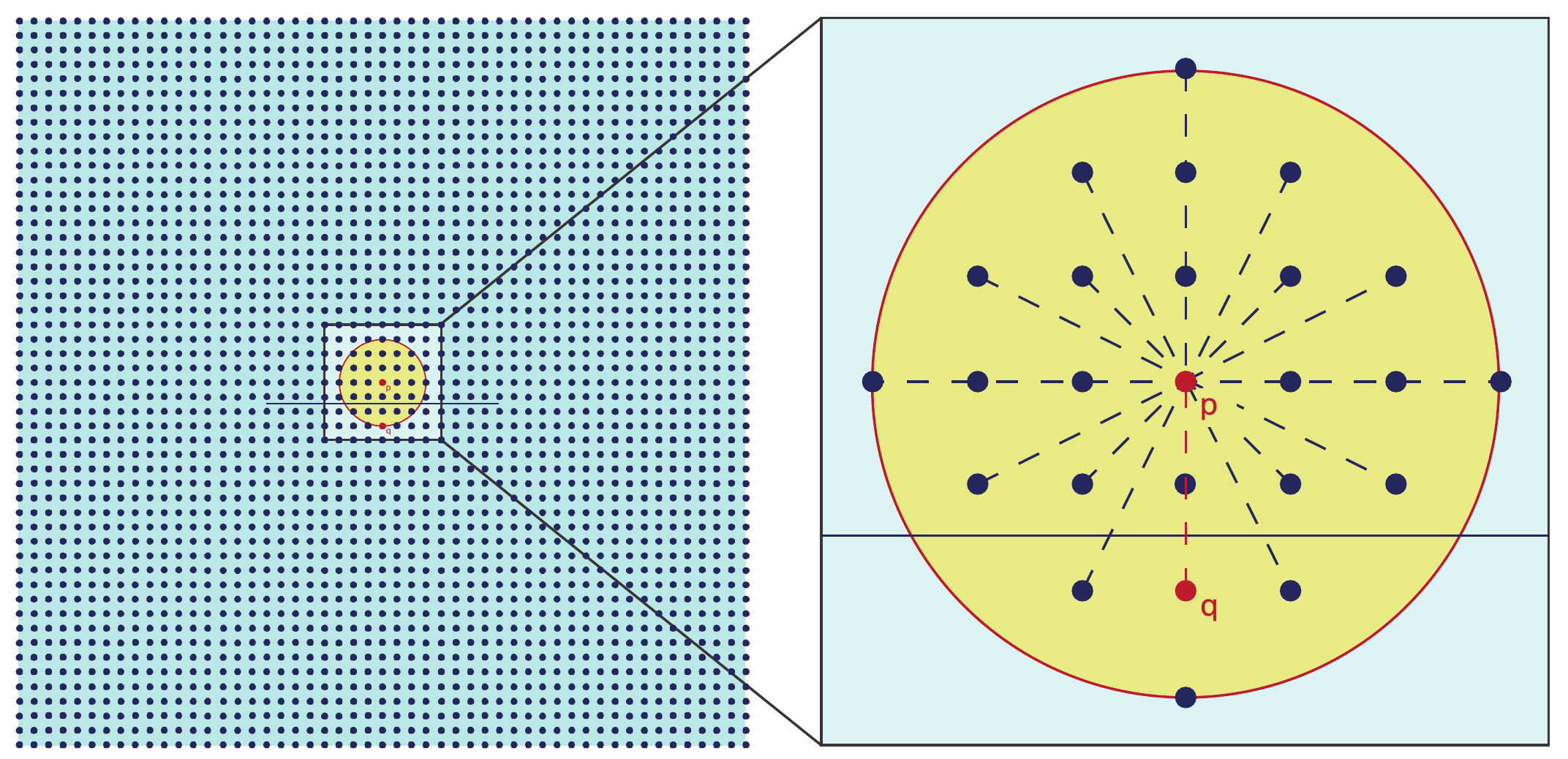}  
	\caption{Schematic diagram of cracks when there is a fracture on bond between $\mathbf{x}_p$ and $\mathbf{x}_q$: the left part is peridynamic discretization, and the right part is geometry details} 
	\label{Fig:4}
\end{figure}
Note that $A_{p,p}$ is the sum of $A_{p,q}$, thus matrix entries $\hat{A}_{p,p}$ are also affected by the broken bonds. For the purpose of explaining $A_{p,p}$, we introduce a set $\mathbf{v}_p$ to store $\mathbf{x}_q$ on the broken bond, which can be expressed as $\mathbf{v}_p=\{\x_q:\text{the bonds between $\mathbf{x}_p$ and $\mathbf{x}_q$ are broken}\}$. The entry $A_{p,p}$ can be given by:
\begin{equation}
\begin{split}
	A_{p,p}=-\sum_{\mathbf{x}_{q}\in\mathcal{B}_{\delta}(\mathbf{x}_p)\cap\Omega}A_{p,q}+\sum_{\mathbf{x}_{q}\in\mathbf{v}_{p}}A_{p,q}
\end{split}
\end{equation}

 Then the matrix $\mathbf{A}$ can be decomposed into a form:
    \begin{equation}
    \mathbf{A}=\mathbf{\hat{A}}+\mathbf{D}^{f},
    \end{equation}
    where

   \begin{equation}\label{df}  
   D_{p,q}^f=  
   \begin{cases}
    -\hat{A}_{p,q}& \text{if $\mathbf{x}_p\in\Omega_{br}$, $\mathbf{x}_q\in\mathbf{v}_p$}
    \\
    \sum\limits_{\mathbf{x}_s\in\mathbf{v}_p} \hat{A}_{p,s}&\text{if $\mathbf{x}_p\in\Omega_{br}$, $q=p$}\\
    0,&\text{otherwise}
    \end{cases}
   \end{equation}
   
   Thus $\mathbf{f}$ can be expressed as:
  \begin{equation}\label{tea:fr}
  f_p=
   \begin{cases}
    f^{'}_p+\sum\limits_{\x_q\in\mathbf{v}_p} \hat{A}_{p,q}(u_p-u_q),&\text{if $\mathbf{x}_p\in\Omega_{br}$, }\\
    f^{'}_p&\text{otherwise}
   \end{cases}
  \end{equation}

  The calculation of form $\mathbf{D}^f\mathbf{u}$ is related to the number of broken bonds.  In fact, the cracks are lower dimensional manifolds compared to the domain’s dimension, which means that the number of material points on $\Omega_{br}$ will not exceed $N^{\frac{d-1}{d}}$. Thus the calculation is generally obtained  by $N^{\frac{d-1}{d}}N$, which is $O(N^{\frac{2d-1}{d}})$.

  After calculating the matrix $\mathbf{D}^f$, we will consider the influence of $\Omega_e$, $\Omega_{in}$, which means that stiff matrix $A$ is decomposed into the form $\mathbf{A=\hat{A}+D^f+D^e}$. If $\Omega_c\cap\Omega_{br}\neq\varnothing$, we also replace the corresponding displacement in the iteration.

\section{A fast matrix-based method for a 3D  linear bond-based peridynamic model}

  To develop the MSBFM on the 3D model, a linear bond-based peridynamics in three spaces dimensions on a block $\Omega=\Omega_c\cap\Omega_s=[x_l,x_r]\times[y_b,y_t]\times[z_c,z_d]$ is introduced in this section. Here $\Omega_c$ represents the area affected by volume constrained boundary conditions. Identical to the form (\ref{fxy}),  $\mathbf{x}_p$, $\mathbf{u}_p$ and $\mathbf{f}_p$ can be denoted as $\mathbf{x}_p=[x_p,y_p,z_p]^T$, $\mathbf{u}_p=[u_p^x,u_p^y,u_p^z]^{T},\mathbf{f}_p=[f_p^x,f_p^y,f_p^z]^{T},1\leq p \leq N$. Then $\mathbf{f}$ obtained by Eq.(\ref{mat:f1}) can be expressed as:
  \begin{equation}\label{thm:e1}
  \left[
  \begin{matrix}
  \mathbf{f}^x\\
  \mathbf{f}^y\\
  \mathbf{f}^z\\
  \end{matrix}
  \right]=
  \left[
  \begin{matrix}
  \mathbf{A}^{xx}&\mathbf{A}^{xy}&\mathbf{A}^{xz}\\
  \mathbf{A}^{xy}&\mathbf{A}^{yy}&\mathbf{A}^{yz}\\
  \mathbf{A}^{xz}&\mathbf{A}^{yz}&\mathbf{A}^{zz}
  \end{matrix}
  \right]
  \left[
  \begin{matrix}
  \mathbf{u}^x\\
  \mathbf{u}^y\\
  \mathbf{u}^z\\
  \end{matrix}
  \right]  	  	
  \end{equation}
  Here $\mathbf{f}^z=[f_1^z,\dots,f_p^z]^T$, $\mathbf{u}^z=[u_1^z,\dots,u_p^z]^T$, and $\mathbf{A}^{xz}$, $\mathbf{A}^{yz}$, $\mathbf{A}^{zz}$ are defined as:
  \begin{equation}\label{mae:e1}
  \begin{split}
    &A_{p,q}^{xx}=
  \begin{cases}
  \alpha\mu \dfrac{(x_{q}-x_{p})^2}{((x_{q}-x_{p})^{2}+(y_{q}-y_{p})^{2}+(z_{q}-z_{p})^{2})^{\frac{3}{2}}}\lambda_{q} V_{q}& q\neq p,\forall \mathbf{x}_p \in \Omega_s
  \\-\sum\limits_{\mathbf{x}_{q}\in\mathcal{B_{\delta}}(\mathbf{x}_p)\cap\Omega}\alpha\mu \dfrac{(x_{q}-x_{p})^2}{((x_{q}-x_{p})^{2}+(y_{q}-y_{p})^{2}+(z_{q}-z_{p})^{2})^{\frac{3}{2}}}\lambda_{q} V_{q}& q=p,\forall \mathbf{x}_p \in \Omega_s
  \end{cases}\\
   &A_{p,q}^{xy}=
 \begin{cases}
 \alpha\mu \dfrac{(x_{q}-x_{p})(y_{q}-y_{p})}{((x_{q}-x_{p})^{2}+(y_{q}-y_{p})^{2}+(z_{q}-z_{p})^{2})^{\frac{3}{2}}}\lambda_{q} V_{q}& q\neq p,\forall \mathbf{x}_p \in \Omega_s
 \\-\sum\limits_{\mathbf{x}_{q}\in\mathcal{B_{\delta}}(\mathbf{x}_p)\cap\Omega}\alpha\mu \dfrac{(x_{q}-x_{p})(y_{q}-y_{p})}{((x_{q}-x_{p})^{2}+(y_{q}-y_{p})^{2}+(z_{q}-z_{p})^{2})^{\frac{3}{2}}}\lambda_{q} V_{q}& q=p,\forall \mathbf{x}_p \in \Omega_s
 \end{cases}\\
  &A_{p,q}^{xz}=
  \begin{cases}
  \alpha\mu \dfrac{(x_{q}-x_{p})(z_{q}-z_{p})}{((x_{q}-x_{p})^{2}+(y_{q}-y_{p})^{2}+(z_{q}-z_{p})^{2})^{\frac{3}{2}}}\lambda_{q} V_{q}& q\neq p,\forall \mathbf{x}_p \in \Omega_s
  \\-\sum\limits_{\mathbf{x}_{q}\in\mathcal{B_{\delta}}(\mathbf{x}_p)\cap\Omega}\alpha\mu \dfrac{(x_{q}-x_{p})(z_{q}-z_{p})}{((x_{q}-x_{p})^{2}+(y_{q}-y_{p})^{2}+(z_{q}-z_{p})^{2})^{\frac{3}{2}}}\lambda_{q} V_{q}& q=p,\forall \mathbf{x}_p \in \Omega_s
  \end{cases}\\
   &A_{p,q}^{yy}=
  \begin{cases}
  \alpha\mu \dfrac{(y_{q}-y_{p})^2}{((x_{q}-x_{p})^{2}+(y_{q}-y_{p})^{2}+(z_{q}-z_{p})^{2})^{\frac{3}{2}}}\lambda_{q} V_{q}& q\neq p,\forall \mathbf{x}_p \in \Omega_s
  \\-\sum\limits_{\mathbf{x}_{q}\in\mathcal{B_{\delta}}(\mathbf{x}_p)\cap\Omega}\alpha\mu \dfrac{(y_{q}-y_{p})^2}{((x_{q}-x_{p})^{2}+(y_{q}-y_{p})^{2}+(z_{q}-z_{p})^{2})^{\frac{3}{2}}}\lambda_{q} V_{q}& q=p,\forall \mathbf{x}_p \in \Omega_s
  \end{cases}\\
  &A_{p,q}^{yz}=
  \begin{cases}
  \alpha\mu \dfrac{(y_{q}-y_{p})(z_{q}-z_{p})}{((x_{q}-x_{p})^{2}+(y_{q}-y_{p})^{2}+(z_{q}-z_{p})^{2})^{\frac{3}{2}}}\lambda_{q} V_{q}& q\neq p,\forall \mathbf{x}_p \in \Omega_s
  \\-\sum\limits_{\mathbf{x}_{q}\in\mathcal{B_{\delta}}(\mathbf{x}_p)\cap\Omega}\alpha\mu \dfrac{(y_{q}-y_{p})(z_{q}-z_{p})}{((x_{q}-x_{p})^{2}+(y_{q}-y_{p})^{2}+(z_{q}-z_{p})^{2})^{\frac{3}{2}}}\lambda_{q} V_{q}& q=p,\forall \mathbf{x}_p \in \Omega_s
  \end{cases}\\
  &A_{p,q}^{zz}=
  \begin{cases}
  \alpha\mu \dfrac{(z_{q}-z_{p})^2}{((x_{q}-x_{p})^{2}+(y_{q}-y_{p})^{2}+(z_{q}-z_{p})^{2})^{\frac{3}{2}}}\lambda_{q} V_{q}& q\neq p,\forall \mathbf{x}_p \in \Omega_s
  \\-\sum\limits_{\mathbf{x}_{q}\in\mathcal{B_{\delta}}(\mathbf{x}_p)\cap\Omega}\alpha\mu \dfrac{(z_{q}-z_{p})^2}{((x_{q}-x_{p})^{2}+(y_{q}-y_{p})^{2}+(z_{q}-z_{p})^{2})^{\frac{3}{2}}}\lambda_{q} V_{q}& q=p,\forall \mathbf{x}_p \in \Omega_s
  \end{cases}\\
  \end{split}
  \end{equation}
 According to the symmetry of the kernel function, we can get $\mathbf{A}^{yx}=\mathbf{A}^{xy}$, $\mathbf{A}^{zx}=\mathbf{A}^{xz}$ and $\mathbf{A}^{zy}=\mathbf{A}^{yz}$.
  Here we only consider $\mathbf{f}^x=\mathbf{A}^{xx}\mathbf{u}^x$, and record it as $\mathbf{f=Au}$.

  This model can also be discretized by uniform mesh, as shown in 2D model. Here the material points can be expressed as $\mathbf{x}_p=[x_i,y_j,z_k]^T$, where $x_i=x_l+(i-1/2)h_x$, $y_j=y_b+(j-1/2)h_y$, $z_k=z_c+(k-1/2)h_z$. $h_x$, $h_y$, $h_z$ are positive constants representing the grid spacing, and we let $h_x=h_y=h_z=h$. $z$ is the index of the layer, and $N_z$ denotes the numbers of  intervals in the $z$ directions. Node number $p$ can be obtained  by $p=(k-1)N_xN_y+(j-1)N_x+i$.  Thus $\mathbf{f}$ and $\mathbf{u}$ can be rewritten as:  
  \begin{equation}
  \begin{split}
  &\mathbf{u}=[u_{1,1,1},\dots,u_{N_x,1,1},\dots,u_{1,N_y,1},\dots,u_{N_x,N_{y},1},\dots, u_{1,1,N_z},\dots,u_{N_x,N_{y},N_{z}}]^{T}\\
  &\mathbf{f}=[f_{1,1,1},\dots,f_{N_x,1,1},\dots,f_{1,N_y,1},\dots,f_{N_x,N_{y},1},\dots, f_{1,1,N_z},\dots,f_{N_x,N_{y},N_{z}}]^{T}
  \end{split}
  \end{equation}    
  where the matrix entry $A_{p,q}$ can be written as:
  \begin{equation}
  A_{p,q}=
  \begin{cases}
  \alpha\mu \dfrac{(x_{i^{\prime}}-x_{i})^{2}}{((x_{i^{\prime}}-x_{i})^{2}+(y_{j^{\prime}}-y_{j})^{2}+(z_{k^{\prime}}-z_{k})^{2})^{\frac{3}{2}}}\lambda_{i^{'},j^{'},k^{'}}V_{i^{'},j^{'},k^{'}}& q\neq p, \forall \mathbf{x}_p\in\Omega_s
  \\-\sum\limits_{\mathbf{x}_{q}\in\mathcal{B}_{\delta}(\mathbf{x}_p)\cap\Omega}\alpha\mu \dfrac{(x_{i^{\prime}}-x_{i})^{2}}{((x_{i^{\prime}}-x_{i})^{2}+(y_{j^{\prime}}-y_{j})^{2}+(z_{k^{\prime}}-z_{k})^{2})^{\frac{3}{2}}}\lambda_{i^{'},j^{'},k^{'}}V_{i^{'},j^{'},k^{'}}& q=p, \forall \mathbf{x}_p\in\Omega_s
  \end{cases}
  \end{equation}
  Here $V_{i^{'},j^{'},k^{'}}=h^3$.

  Following the spatial partition in the 2D model, we divide the region  $\Omega_s$ into $\Omega_{in}$, $\Omega_e$ and $\Omega_{br}$. They represent the internal area, the area with the incomplete disk, and the area with broken bonds, respectively.

  Then $\mathbf{A}$ becomes a stiff matrix  with a block-banded-block-banded-block structure, which means $\mathbf{A}$ is a matrix composed of $N_z^2$ matrix blocks $\mathbf{\bar{B}}_{i,j}$. Each matrix block $\mathbf{\bar{B}}_{i,j}$ represents the action of the $j-$th layer  on the $i-$th layer , and the structure of $\mathbf{\bar{B}}_{i,j}$ is a block-banded-block matrix, which is as same as the form (\ref{stu:bad}), see Fig. \ref{fig:5}.

    \begin{figure}[h]\
  	\centering            
  	\includegraphics[scale=0.6]{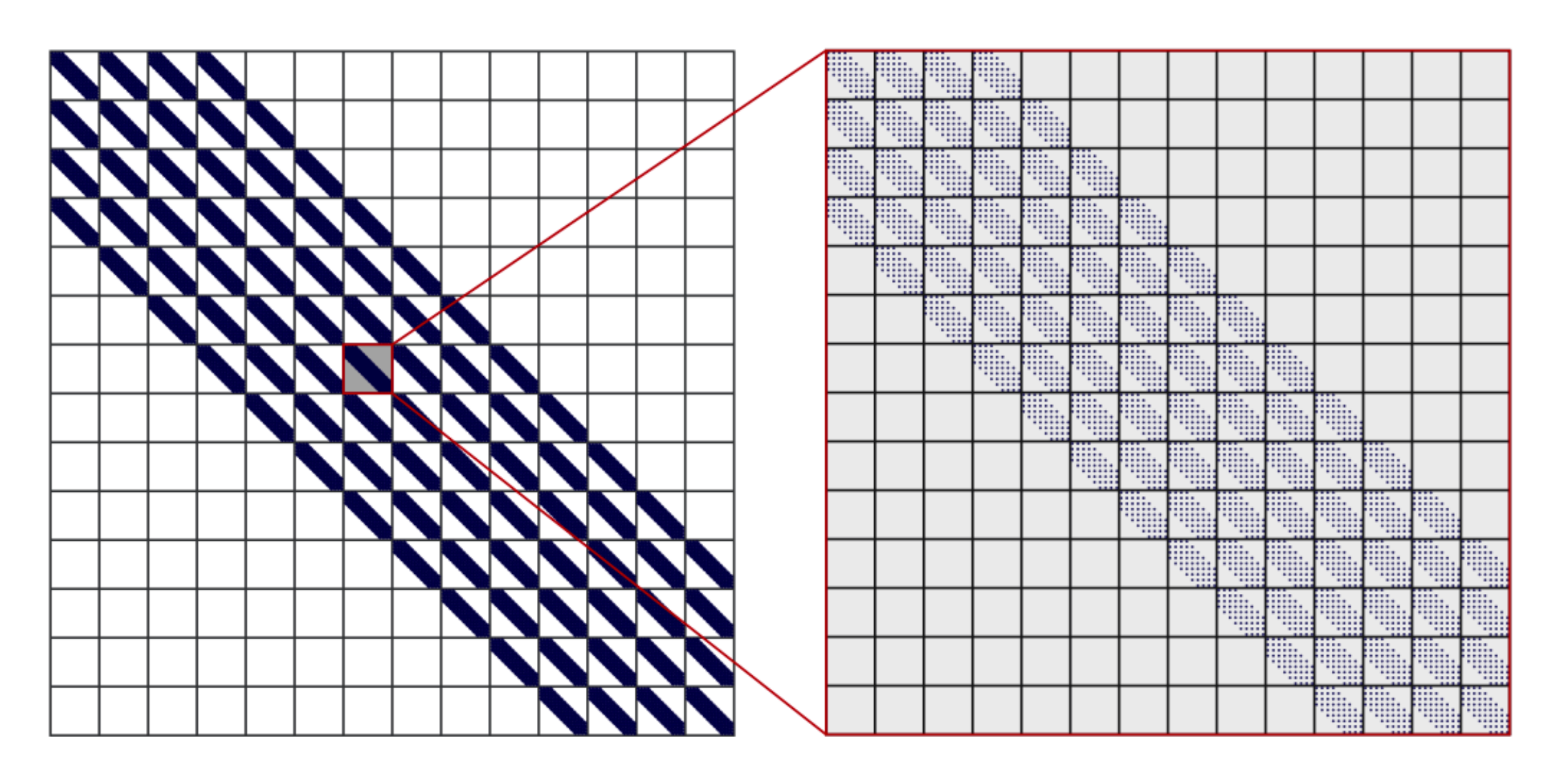}  
  	\caption{Illustration of the matrix structure in 3D model:The left side is the matrix $\mathbf{A}$ composed of matrix block $\mathbf{\bar{B}}_{i,j}$, and the right side is the structure of matrix block $\mathbf{\bar{B}}_{i,j}$} 
  	\label{fig:5}
  \end{figure}
  
  Similar to Eq.(\ref{frb:e1}) and Eq.(\ref{frb:e2}), $A_{p,q}$ can be proved to satisfy the following properties without considering broken bonds:
  \begin{equation}\label{thm:e2}
  \begin{split}
  A_{p,q}&=\alpha\frac{(x_{i^{\prime}}-x_{i})^{2}}{((x_{i^{\prime}}-x_{i})^{2}+(y_{j^{\prime}}-y_{j})^{2}+(z_{k^{\prime}}-z_{k})^{2})^{\frac{3}{2}}}\lambda_{i^{'},j^{'}}h^3\\
  &=\alpha\frac{(i^{\prime}-i)^{2}}{((i^{\prime}-i)^{2}+(j^{\prime}-j)^{2}+(k^{\prime}-k)^{2})^{\frac{3}{2}}}\lambda_{i^{'},j^{'}}h^2\\ 
  \end{split} 
  \end{equation}
 where  $\mathbf{x}_{q}=[x_{i^{'}},y_{j^{'}},z_{k^{'}}]^T$, $q\neq p$. When $q=p$, we can also get:
  \begin{equation}\label{thm:e3}
  \begin{split}
  A_{p,p}&=-\sum\limits_{\mathbf{x}_{q}\in\mathcal{B}_{\delta}(\mathbf{x}_p)\cap\Omega} \alpha\frac{(x_{i^{\prime}}-x_{i})^{2}}{((x_{i^{\prime}}-x_{i})^{2}+(y_{j^{\prime}}-y_{j})^{2}+(z_{k^{\prime}}-z_{k})^{2})^{\frac{3}{2}}}\lambda_{i^{'},j^{'}}h^3\\&=-\sum\limits_{\mathbf{x}_{q}\in\mathcal{B}_{\delta}(\mathbf{x}_p)\cap\Omega}\alpha\frac{(i^{\prime}-i)^{2}}{((i^{\prime}-i)^{2}+(j^{\prime}-j)^{2}+(k^{\prime}-k)^{2})^{\frac{3}{2}}}\lambda_{i^{'},j^{'}}h^2
  \end{split}
  \end{equation}
  Identical to the 2D model, we can prove that matrix $\mathbf{A}$  is only related to $i^{'}-i$, $j^{'}-j$, $k^{'}-k$ and introduce a $2M+1$-by-$2M+1$-by-$2M+1$ tensor $\mathbf{K}$ to store entries $A_{p,q}$, namely:
  \begin{equation}
  A_{p,q}=K_{m,n,l},\quad 0\le m,n,l \le 2M+1, \quad \forall \mathbf{x}_p=(x_i,y_j,z_k)\in\Omega_{in}
  \end{equation}
  Here $m=i^{'}-i+M+1$, $n=j^{'}-j+M+1$, $l=k^{'}-k+M+1$.
   Then a Toeplitz-Block-Toeplitz-Block-Toeplitz(TBTBT) matrix $\mathbf{\hat{A}}$  can also be defined as:
  \begin{equation}
  \hat{A}_{p,q}= K_{m,n,l}, \quad \forall \mathbf{x}_p\in\Omega
  \end{equation}
  
 The form $\mathbf{f^{'}}=\mathbf{\hat{A}u}$  can be solved by a fast tensor-tensor multiplication(FTTM), as shown in Algorithm \ref{vv4}:
  \begin{algorithm}   	
 	\SetAlgoNoLine 
 	\SetKwInOut{Input}{\textbf{Input}}\SetKwInOut{Output}{\textbf{Output}}   	
 	$\mathbf{\hat{G}}=\mathbf{FFT3}(\mathbf{G}),\mathbf{\hat{U}}=\mathbf{FFT3}(\mathbf{U}$) \\	\vspace{0.03in}
 	$\mathbf{\hat{H}}=\mathbf{\hat{G}}\circ\mathbf{\hat{U}}$ \\ \vspace{0.03in}
 	$\mathbf{H}=\mathbf{FFT3^{-1}}(\mathbf{\hat{H}})$\\\vspace{0.03in}
 	\caption{\textbf{FMVM}\label{vv4}}
 \end{algorithm}

   Here we use $\mathbf{FFT3}$ and $\mathbf{FFT3^{-1}}$ to denote three-dimension FFT and iFFT operations. $\mathbf{H}$ is a $2N_x$-by-$2N_y$-by-$2N_z$ tensor and we can obtain $\mathbf{f}$ by $f_{kN_xN_y+jN_x+i}=H_{i,j,k}$, if $1\leq i\leq N_x$,  $1\leq j\leq N_y$, $1\leq k\leq N_z$. For the first column $\mathbf{G}$ of the extended matrix embedded by the tensor $\mathbf{K}$, we define it as
a $2N_x$-by-$2N_y$-by-$2N_z$ tensor, namely:
 \begin{equation}
 G_{i,j,k}=
 \begin{cases}
 0&i\in[M+2,2N_x-M],j\in[M+2,2N_y-M],k\in[M+2,2N_z-M]\\
 K_{m,n,l}&\text{otherwise}\\
 \end{cases}  	
 \end{equation}
 where 
 \begin{equation}
 \begin{split}
 &m=\begin{cases}
 i+M&\text{$i\in[1,M+1]$}\\
 i-2N_x+M&\text{$i\in[2N_x-M+1,2N_x]$}
 \end{cases}\\
 &n=\begin{cases}
 j+M&\text{$j\in[1,M+1]$}\\
 j-2N_y+M&\text{$j\in[2N_y-M+1,2N_y]$}
 \end{cases}\\
 &l=\begin{cases}
 k+M&\text{$k\in[1,M+1]$}\\
 k-2N_z+M&\text{$k\in[2N_z-M+1,2N_z]$}
 \end{cases}
 \end{split}
 \end{equation}
 The expansion vector $\mathbf{U}$ is also a $2N_x$-by-$2N_y$-by-$2N_z$ tensor, which can be expressed as:
 \begin{equation}
 U_{i,j,k}=\begin{cases}
 u_{(i+1)/2,j,k},&\text{if $i\leq N_x$, $j \leq N_y$, $k\leq N_z$}\\ 
 0,&\text{otherwise}
 \end{cases}
 \end{equation}
 	The codes of the three-dimensional Fourier transform and the inverse Fourier transform in Matlab are called as $\mathbf{\hat{G}}=fftn(\mathbf{G})$ and $\mathbf{{G}}=ifftn(\mathbf{\hat{G}})$.

For the material points in $\Omega_{c}$, $\Omega_{e}$, and $\Omega_{br}$, they are treated in the same way as the matrix mentioned in the form (\ref{stu:bad}), which means the matrix $\mathbf{A}$ of the 3D model is decomposed into $\mathbf{A=\hat{A}+D^e+D^f}$ and the displacement $\mathbf{u}_p$ is replaced with $h(\x_p)$ in each time iteration if $\mathbf{x}_p\in\Omega_{c}$.
 
 The above analysis shows that MSBFM is an algorithm based on the structure of matrix $\mathbf{A}$. Most entries in matrix $\mathbf{A}$ satisfy the TBT or TBTBT structure; thus, the form $\mathbf{f=Au}$ can be accelerated by FFT. For problems in 2D and 3D, broken bonds, incomplete disks, and volume constrained boundary conditions break this structure and need special steps to deal with them.  However, the calculation is $O(N\log N)$ because most material points are in $\Omega_{in}$. Here a flowchart is introduced to illustrate these steps as shown in Fig. \ref{Fig:5}.
  \begin{figure}[h]
	\centering            
	\includegraphics[scale=0.48]{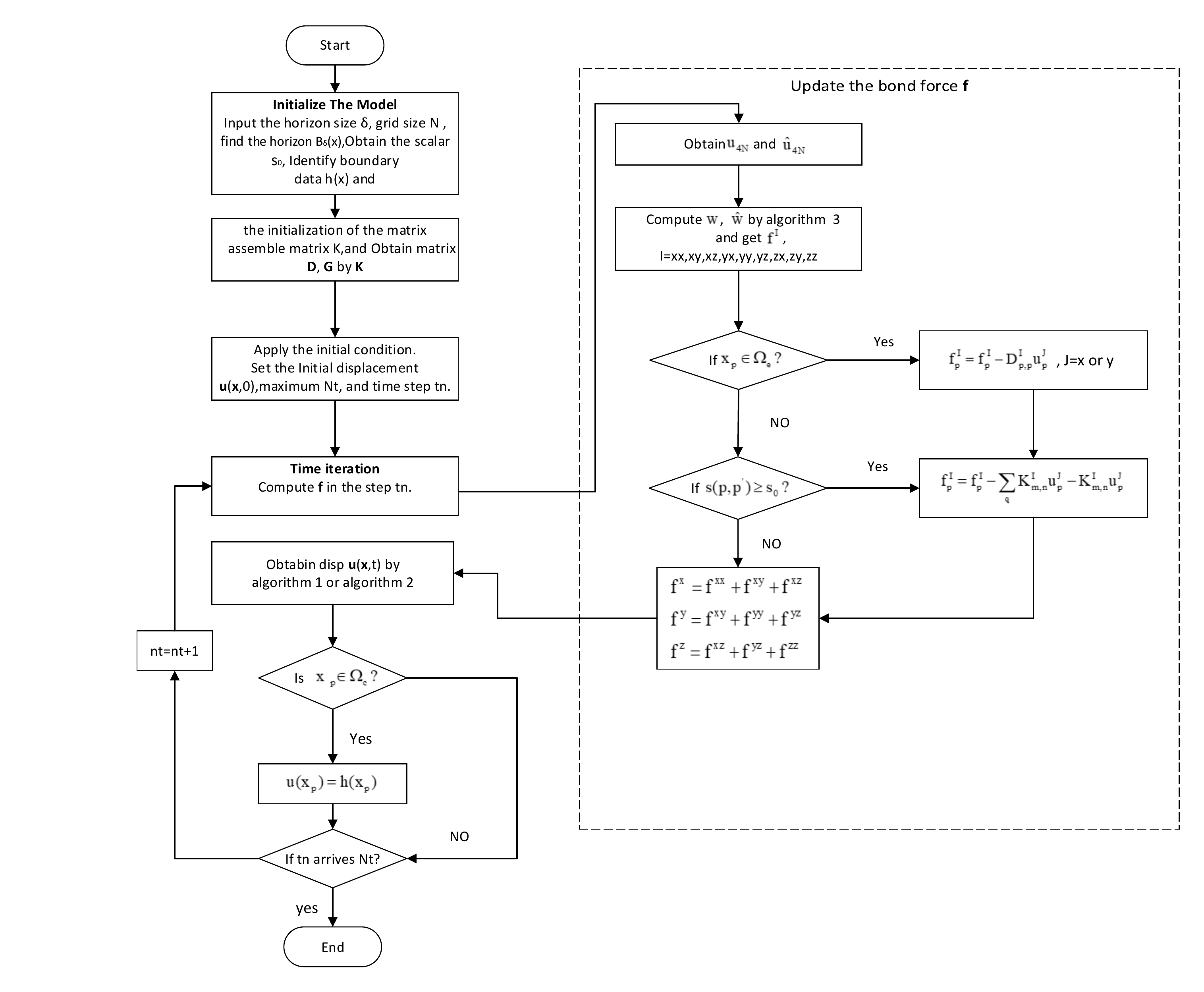}  
	\caption{The flowchart of MSBFM} 
	\label{Fig:5}
\end{figure}

\section{Numerical results}
In this section, the MSBFM is verified by comparing the meshfree method on four examples built on the peridynamic model, including 2D/3D models with various boundary conditions and cracks. We implement these methods in Matlab and run all experiments on a workstation with Intel Xeon Gold 6240(2.6GHz/18C) logical processors and 2048G installed memory. 
 \subsection{Peridynamic in 2D body with external loading}
As a first illustrative example, let us consider a 2D peridynamic model on a plate with external loading in this section.
    
   As shown in Fig. (\ref{Fig:6}), the plate has the width  $W=0.5\mathrm{~m}$, the length  $L=1.0\mathrm{~m}$, and thickness $h=0.0025\mathrm{~m}$. 
        \begin{figure}[h]        
    	\centering  
    	\subfigure[]  	
    	{
    		\begin{minipage}{7cm}
    			\centering          
    			\includegraphics[scale=0.25]{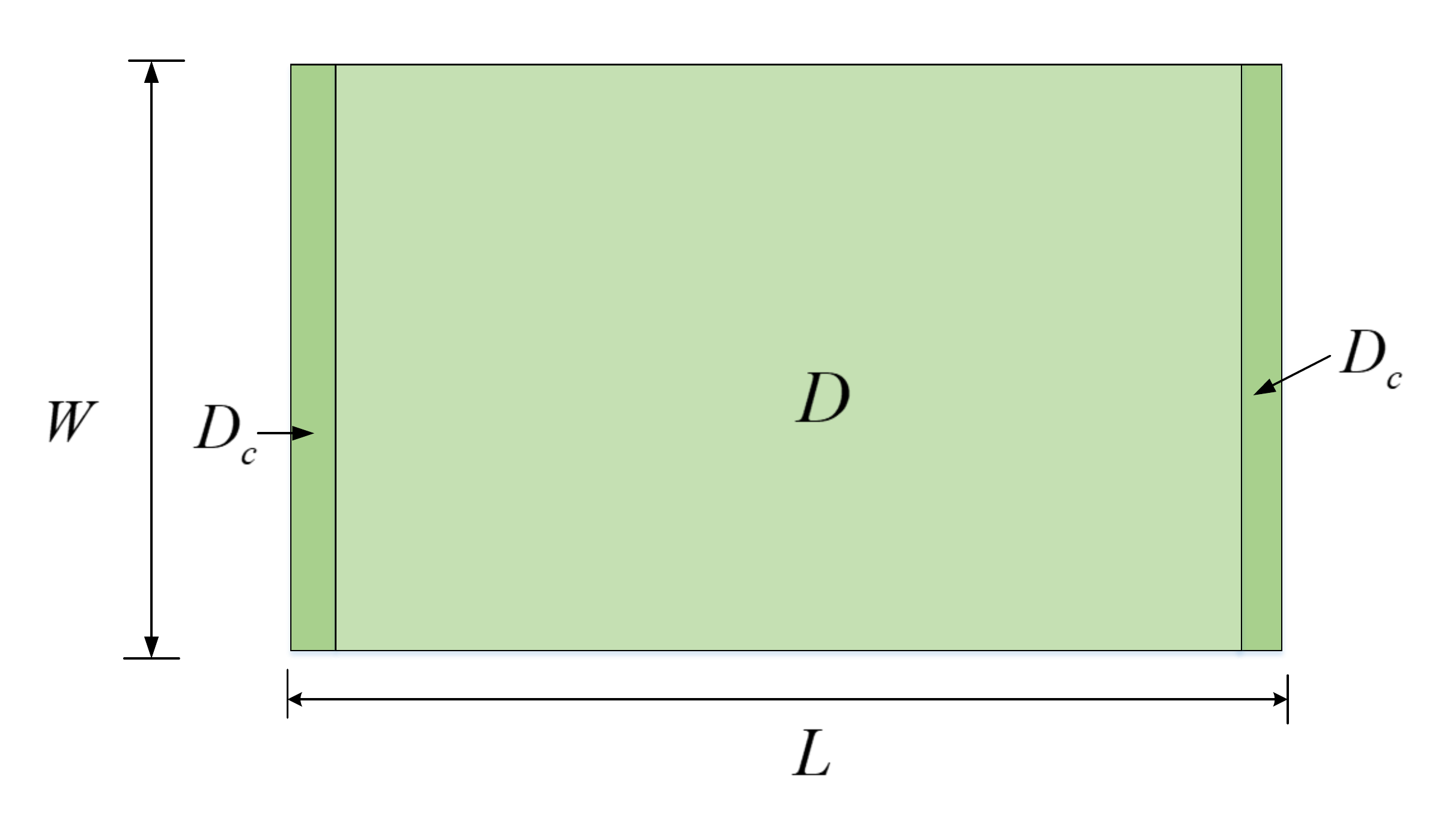}  
    		\end{minipage}
    	}
      \subfigure[]	
    	{
    		\begin{minipage}{8cm}
    			\centering      
    			\includegraphics[scale=0.4]{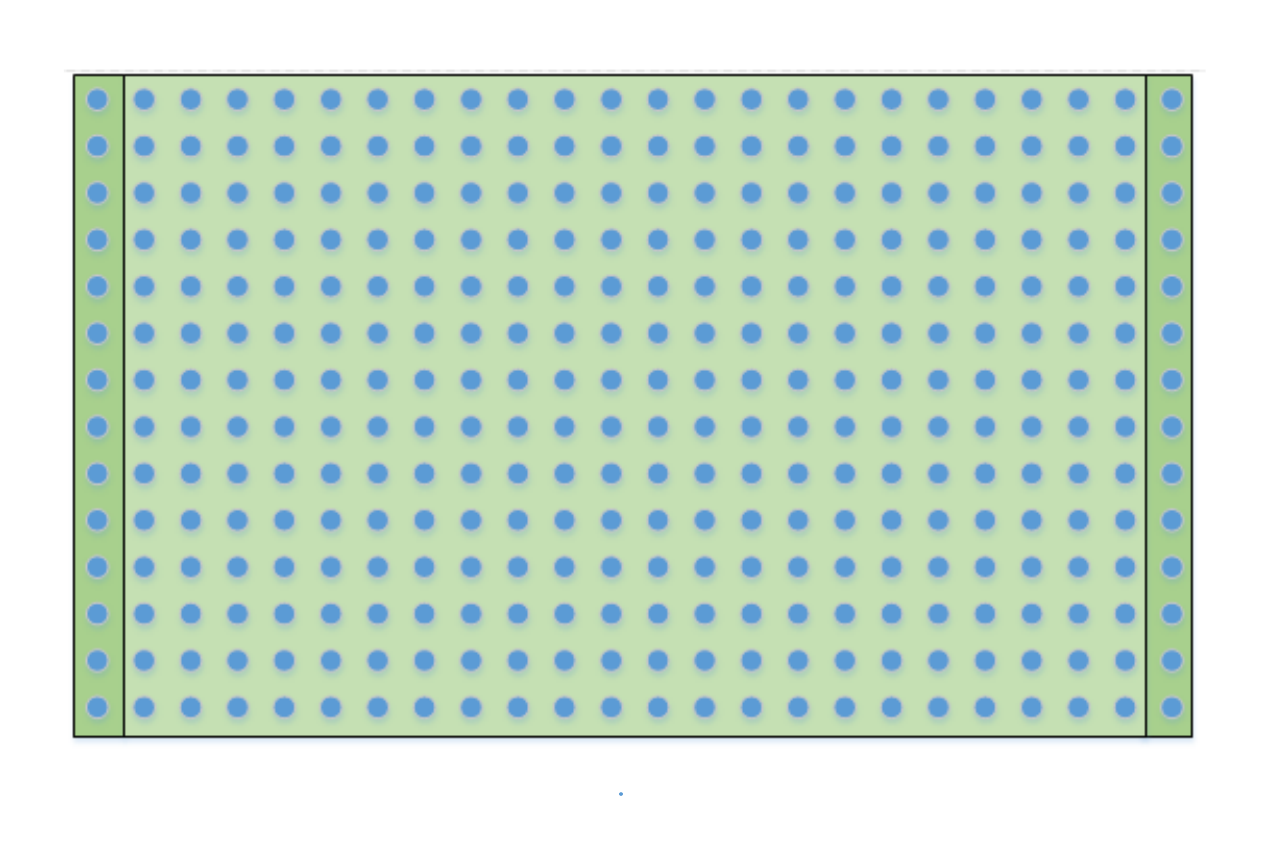}   
    		\end{minipage}
    	}
    	
    	\caption{Geometry of a plate with exteral loads and its discretization (a)plate with external loading; (b)a simple example of uniform grid structure.} 
    	\label{Fig:6}
    \end{figure}
   The material properties are $E=2\times 10^5$Mpa( elastic modulus), $\nu=1/3$(Poisson's ratio) and $\rho=7850kg/\mathrm{~m}^3$(density ). The external loading $b_p=p_0W/h$ is applied to the boundary layer $\mathcal{D}_c$, where the uniaxial tension loading $p_0$ is chosen as $200Mpa$, and the width of $\mathcal{D}_c$ is $h$. $\mathcal{D}$

    We selected the horizon size $\delta$ as $0.03\mathrm{~m}$ and discretize the model by considering $N_x=400$ and $N_y=200$. We can consider implementing MSBFM and the meshfree method in this model since it can be written as a matrix-vector multiplication due to the Eq.(\ref{mmv:fm}) and the ADR method is used for temporal integration(See Section 2.1).

    Fig. \ref{Fig:7}(a) shows the displacement variations obtained by two algorithms when the total time step equals $3000$. It is noticed from Fig. \ref{Fig:7}(a) that the displacement variations by our algorithm has a good match with the results obtained by
    meshfree method when $400\times 200$ mesh elements were employed. 
            \begin{figure}[h]
            	\subfigure[]  	
    	{
    		\begin{minipage}{8cm}  
    			\centering      
    			\includegraphics[scale=0.3]{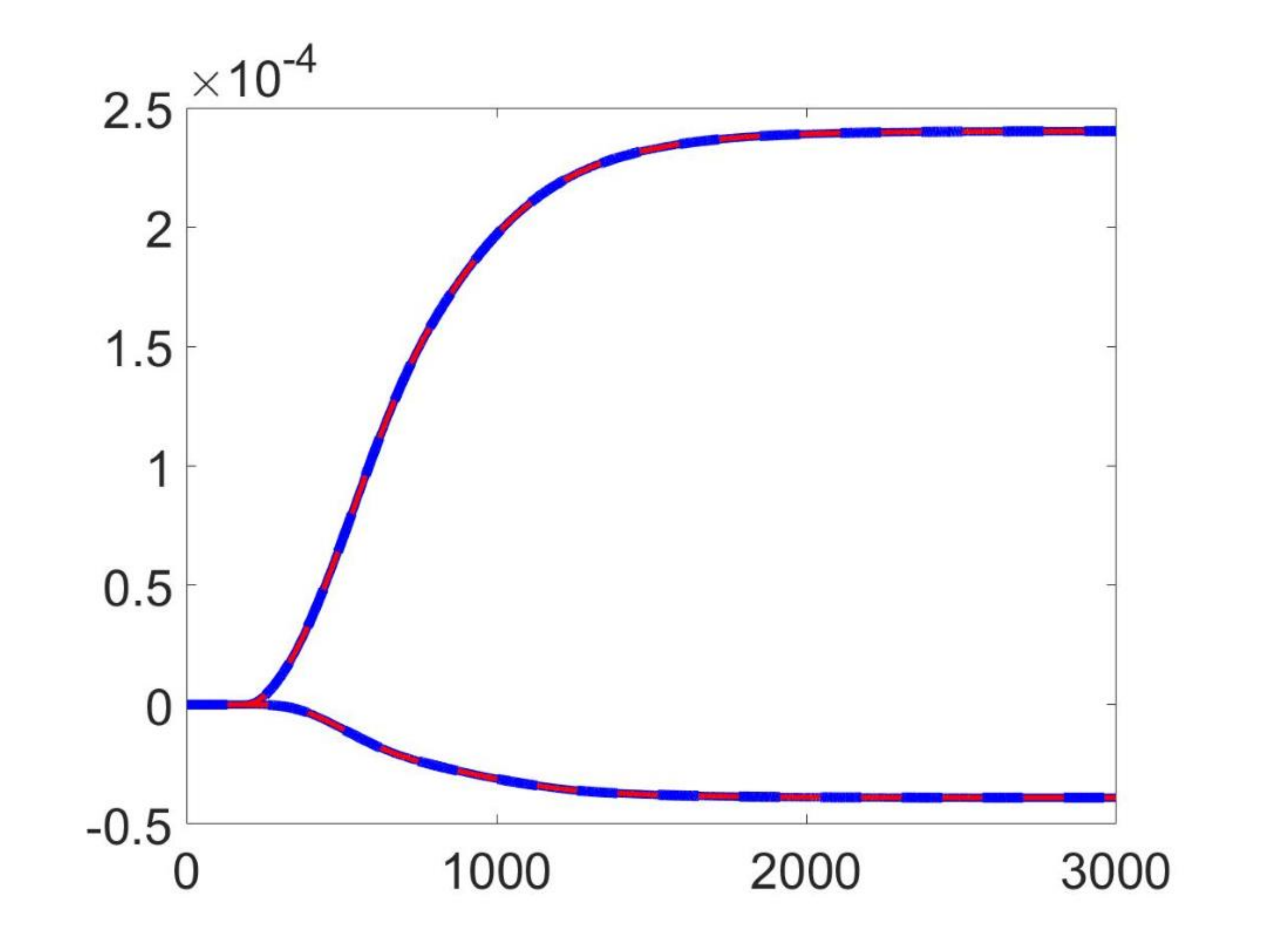}  
    		\end{minipage}
    	}
   \subfigure[]
    	{
    		\begin{minipage}{8cm} 
    			\centering     
    			\includegraphics[scale=0.3]{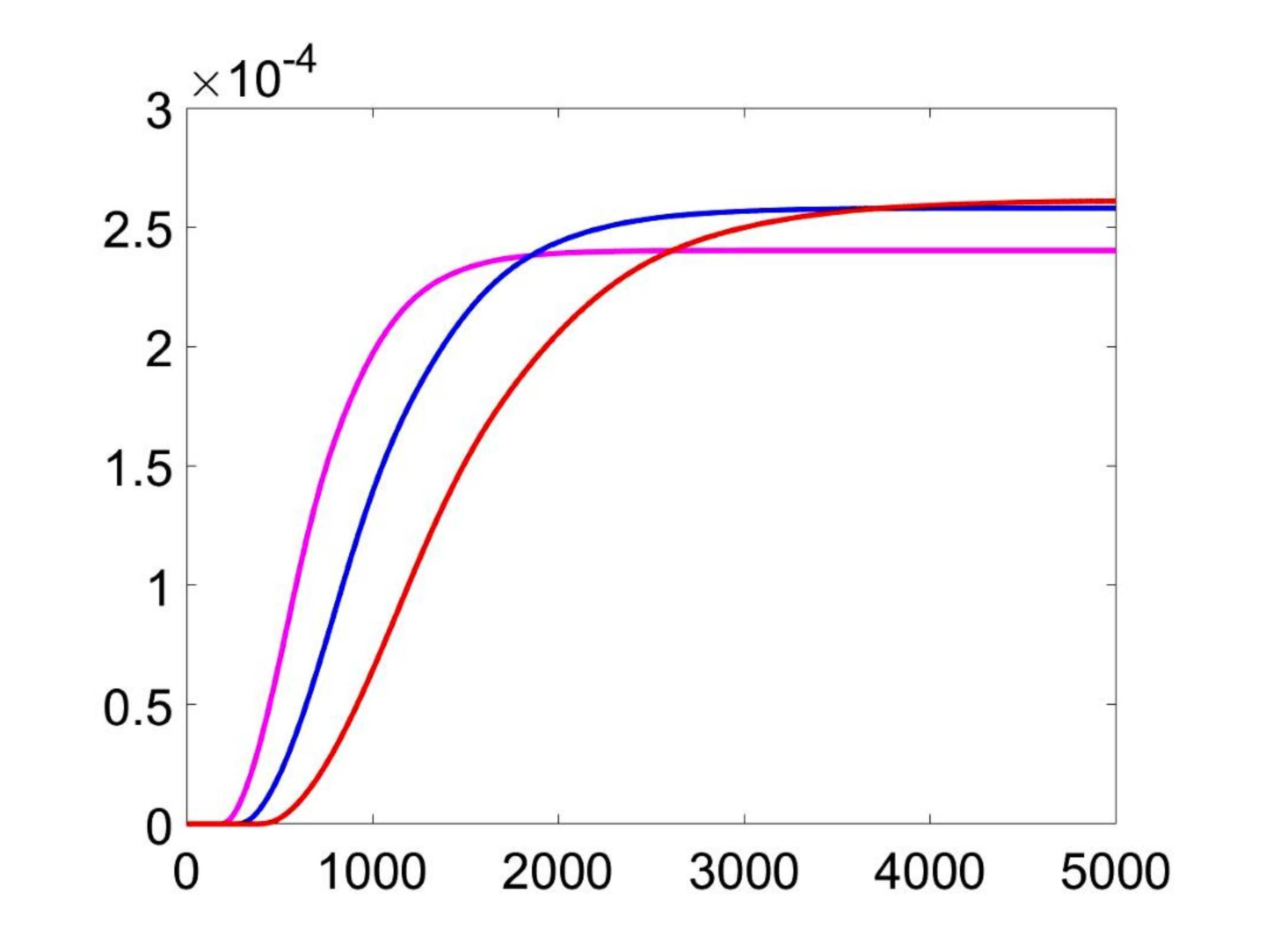}   
    		\end{minipage}
    	}
    	    	\caption{Displacement variation with respect to step number in $x$ direction and $y$ direction for $\mathbf{x}_p=(0.255\mathrm{~m},0.125\mathrm{~m})$: (a) The displacement variations obtained by MSBFM and meshfree method under the mesh $400\times 200$, where the blue solid line represents the data obtained by FMEM, and the blue dotted line represents the data obtained by the meshfree method; (b)comparison of displacements in the $x$-direction under grids is $400\times 200$(green), $800\times 400$(blue), $1600\times 800$(red) by MSBFM.
    	} 
    	\label{Fig:7}
    \end{figure}

\begin{table}\begin{center}
		
		\caption{Performance of meshfree method and MSBFM  in 2D model}
				\vspace{0.15in}
		\setlength{\tabcolsep}{2mm}{
		\begin{tabular}{ccccccc} \hline
			\multicolumn{2}{c}{Mesh} & $400\times 200$ & $800\times400$ & $1600\times 800$ &$3200\times 1600$&$6400\times 3200$\\ \hline
			\multicolumn{2}{c}{Time steps}  & $3000$ & $4000$ & $5000$ &$6000$&$7000$\\ \hline
			
		\multirow{2}*{Meshfree}  
			&Matrix assembly&1m1s&16m40s&4h38m&3d3h&- \\
			&Time stepping&4m35s&4h28m&4d18h&-&-\\\hline
			\multirow{2}*{MSBFM}
			&Matrix assembly&24s&3m9s&40m5s&14h40m&6d6h\\
			&Time stepping &2m30s&11m42s&52m29s&3h40m&1d4h\\ \hline
		\end{tabular}}\label{tab:1}
\end{center}\end{table}
      
    To perform the simulations by using various discretization sizes, we gradually increase the number of grids from $400\times200$ to $1600\times 800$. Table \ref{tab:1} compares the computational time required to perform the simulations using MSBFM and meshfree method . We automatically stop a numerical run if it takes more than 10 days of CPU time. Even in various mesh elements, the results obtained by the two algorithms are the same, so we can only consider one of them when analyzing the properties of material.

    In fact, as the mesh elements increases, the number of time steps required by the ADR method to achieve stability also increases. As shown in the Fig. \ref{Fig:7}(b), when the mesh number increases from $400\times200$ to $1600\times800$, the number of time steps to achieve stability also increases from $1000$ to $3000$. Therefore, we select various time steps according to various mesh elements.

   For the non fracture problems, the calculation is mainly divided into two parts: Phase I: Matrix assembly. For the meshfree method, we need to traverse the horizon of all material points and initialize an $N$-by-$N$ matrix $\mathbf{A}$ in this phase. In the MSBFM method, matrices $\mathbf{K}$ and $\mathbf{D}^e$ are used instead of matrices $\mathbf{A}$, thus reducing the computational time. Phase II: Time stepping. The form $\mathbf{Au}$ and ADR methods are calculated in this phase for meshfree method, while the MSBFM method uses FMVM mentioned in (\ref{vv3}) and form $\mathbf{D^eu}$ to replace the calculation of $\mathbf{Au}$.

    For the part I, the advantages of MSBFM are mainly reflected in two aspects: one is the traversal of horizons of the material points. In this example, we use the method of traversing all material points and comparing distances to find points within the horizon, so the calculation amount is usually $O(N)$. For the MSBFM method, we only need to traverse the information of the material point affected by the incomplete horizon and a complete horizon material point, so we can reduce the calculation to $O(N^{\frac{3}{2}})$. The other is the assembly of stiffness matrix. For the meshfree method, the assembly of stiffness matrix $\mathbf{A}$ need to be considered before the calculation of time integral equation, and the computational complexity of this part is usually $O(N^2)$. The MSBFM method replaces the stiffness matrix $\mathbf{A}$ with matrices $\mathbf{K}$ and $\mathbf{D^e}$, and the calculation can be obtained by $(2M+1)^2=O(N)$ for matrix $\mathbf{K}$, where the calculation of matrix $\mathbf{D}^e$ can be obtained by $N(4\delta/h)=O(N^{\frac{3}{2}})$, neither of which will exceed $O(N^2)$. Table \ref{tab:1} illustrates the time cost for the matrix assembly in meshfree method and MSBFM.
    
     The computational time of the phase II depends on the complexity of the form $\mathbf{f=Au}$, which is $O(N^2)$ and $O(N\log N)$ for meshfree methods and MSBFM. When the number of grids increases from $N_x,N_y$ to $\beta N_x, \beta N_y$, the mesh free time will increase by $\beta^4$ times since $N^2=(\beta^2 N_xN_y)^2=\beta^4 N$. However, the MSBFM time only increases by $\beta^4$ times because $\beta^2N_xN_y\log(\beta^2 N_xN_y)=\beta^2=\beta^2N\log N+2\log\beta$.

  However, the calculation is not strict $O(N\log N)$ for the part II in MSBFM algorithm if we use the surface correction algorithm. In the surface correction algorithm, we need to recalculate $\mathbf{f}_p$ for $\mathbf{x}_p\in\Omega_e$, which is:
  \begin{equation}\label{scr}
  	\mathbf{f}_p=\sum_{\mathbf{x}_q\in\mathcal{B}_{\delta}(\mathbf{x}_p)\cap\Omega}v_{p,q}A_{p,q}(\mathbf{u}_q-\mathbf{u}_p)
  \end{equation}
This means that we need to recalculate the form $\mathbf{Au}$ for the point on $\Omega_e$ by the form (\ref{scr}). In the meshfree method, this part of the calculation is not omitted, but in the MSBFM method, the form $\hat{\mathbf{A}}\mathbf{u}$ and $\mathbf{D^eu}$ is used instead of the original form of $\mathbf{Au}$, thus it will bring an additional calculation, which is $O(N^2)$, for the MSBFM method. In order to store the material points that need to be affected by surface correction, the storage amount will increase to $O(N^2)$, which means that the time of matrix assembly will also increase.

In most cases, the material points on $\Omega_e$ only account for a part of the total material points. Thus the simulation speed of our algorithm is significantly faster, see Table \ref{tab:1.5}.
\begin{table}\begin{center}
		\caption{Performance of meshfree method and MSBFM with surface correction algorithm in 2D model}
				\vspace{0.15in}
		\setlength{\tabcolsep}{2mm}{
			\begin{tabular}{ccccccc} \hline
					\multicolumn{2}{c}{Mesh} & $400\times 200$ & $800\times400$ & $1600\times 800$ &$3200\times 1600$&$6400\times 3200$\\ \hline
				\multicolumn{2}{c}{Time steps}& $3000$ & $4000$ & $5000$ &$6000$&$7000$\\ \hline
				
				\multirow{2}*{Meshfree} 
			    &Matrix assembly&1m1s&16m40s&4h38m&3d3h&- \\
				&Time stepping&4m35s&4h28m&4d18h&-&-\\\hline
				\multirow{2}*{MSBFM} 
				&Matrix assembly&24s&3m9s&40m5s&14h40m&6d6h\\
				&Time stepping &2m40s&22m50s&6h27m&12h48m&-\\ \hline
		\end{tabular}}\label{tab:1.5}
\end{center}\end{table}

 \subsection{Peridynamic in 2D body with a pre-existing crack}
      \begin{figure}[h]
   	\centering  
   	 \subfigure[]  	
   	{
   		\begin{minipage}{6cm}
   			\centering          
   			\includegraphics[scale=0.35]{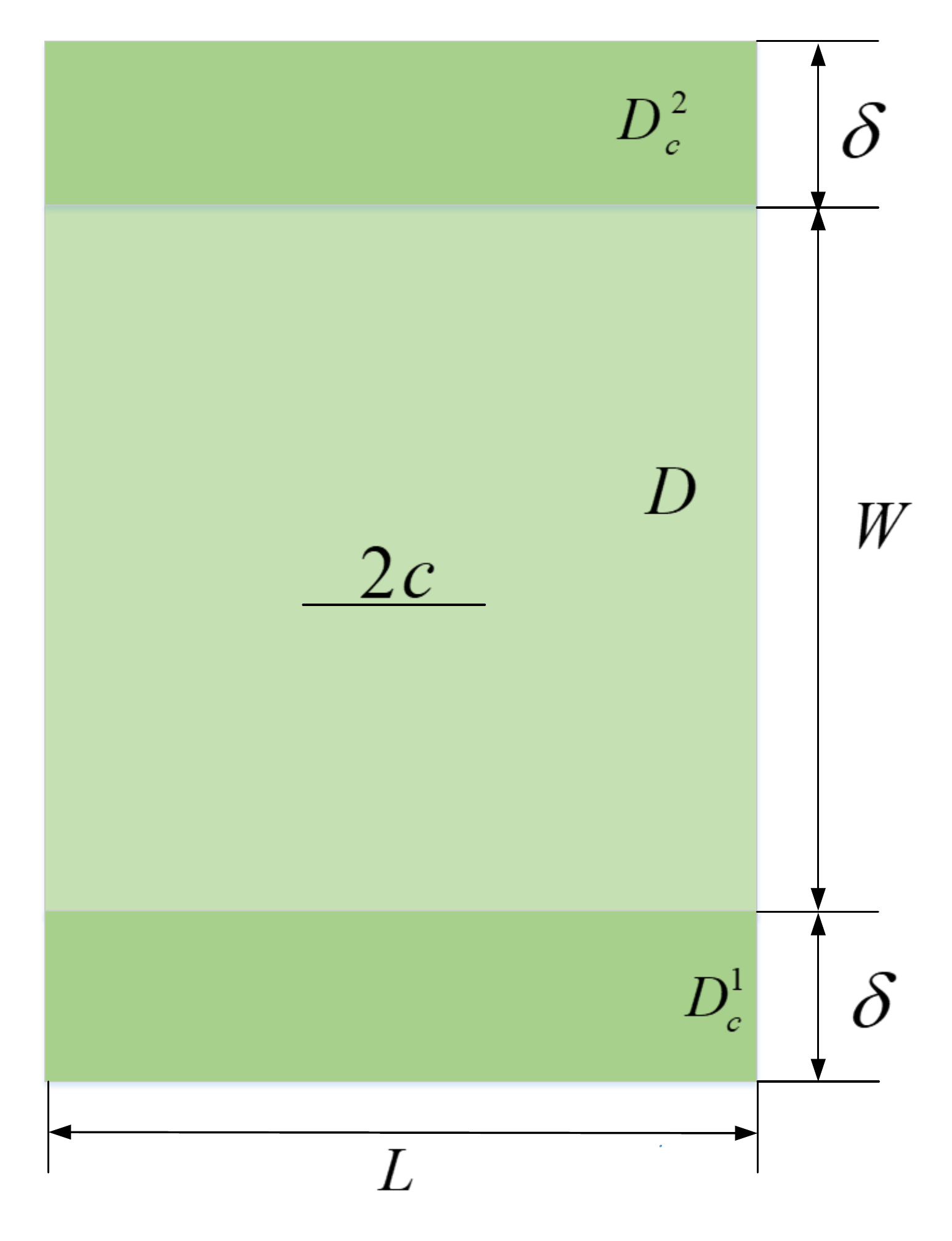}  
   		\end{minipage}
   	}
    \subfigure[]	
   	{
   		\begin{minipage}{6cm}
   			\centering      
   			\includegraphics[scale=0.36]{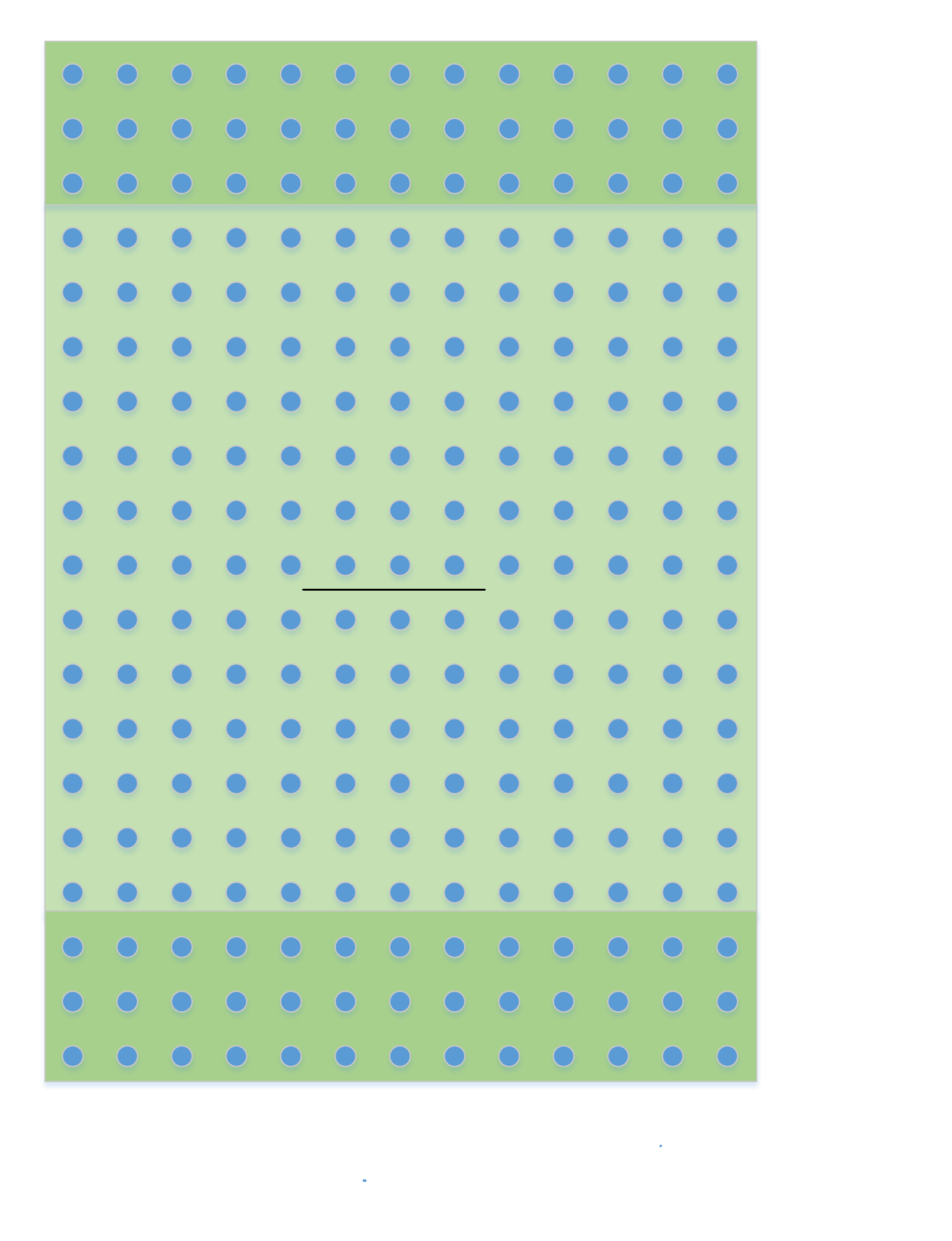}   
   		\end{minipage}
   	}
   	
   	\caption{Geometry of a 2D model with pre-existing crack under velocity constraints and its discretization:(a)plate with pre-existing crack; (b) a simple example of uniform grid structure.} 
   	\label{Fig:8}
   \end{figure}
   
   A 2D model with pre-existing crack is considered in Fig. \ref{Fig:8}. The length $L$ of this plate is  $0.05\mathrm{~m}$,  the width $W$ is $0.05\mathrm{~m}$, and the crack is created as the length $2c=0.01\mathrm{~m}$. Destiny $\rho$, horizon size $\delta$, elastic modulus $E$ and Poisson ratio $\nu$ are chosen to be consistent with Section 5.1.

   Since the PD equation of motion  do not contain any spatial derivatives, the constraints often do not affect the solutions of the integro-differential equations. However, the constraint conditions can still be imposed by introducing a virtual boundary layer, and the displacement on this virtual boundary layer will not be affected by the material points on the actual material area. In this example, virtual boundary layer $\mathcal{D}_c^1$ and $\mathcal{D}_c^2$ with depth $\delta$ is introduced at the upper and lower ends of the actual material area $\mathcal{D}$, and the velocity constraints is applied on $\mathcal{D}_c^1$ and $\mathcal{D}_c^2$, which can be expressed as follows:
   \begin{equation}
   \begin{split}
   	&\dot{u}_y(\mathbf{x}_p,t)=20.0 \mathrm{~m} / \mathrm{s} \quad \mathbf{x}_p\in \mathcal{D}_c^1\\
   	&\dot{u}_y(\mathbf{x}_p,t)=-20.0 \mathrm{~m} / \mathrm{s} \quad \mathbf{x}_p\in \mathcal{D}_c^2
   \end{split}   
   \end{equation}
   
   Although the displacement of the material point on the virtual material layer is independent of the actual problem, we still consider its displacement to ensure that the stiffness matrix $\mathbf{A}$ calculated in the MSBFM algorithm is a square matrix. For the material point $\x_p$ on $\mathcal{D}_c^1$ and $\mathcal{D}_c^2$, displacement $\mathbf{u}_p$ is calculated from two aspects: (a) In the $x$ direction, $\x_p$ is treated as the material point on $\Omega_c$, which means that we need to replace the displacement $u_p^x$ after computing $\mathbf{f^x=A^{xx}u^x+A^{xy}u^x}$. (b) In the $y$ direction, we treat $\x_p$ as the material point on $\Omega_e$, which means that we need to subtract the corresponding $\mathbf{D}^e$ mentioned in (\ref{mat:se}) when calculating $\mathbf{f}^{xx}=\mathbf{A}^{yx}\mathbf{u}^{y}+\mathbf{A}^{yy}\mathbf{u}^{y}$.

We discretize the model with a grid size of $600\times600$, and the Velocity Verlet algorithm is chosen for time discretization because this is a time-dependent problem. We collect data from from $t=0$ to $t=1350$ with a time-step size of $\Delta t = 1.3367\times 10^{-8}s$. Fig. \ref{Fig:9} shows the crack simulation under two algorithms.
  \begin{figure}[h]
  	\centering    	
  	\subfigure[]
  	{
  		\begin{minipage}{7.5cm}
  			\centering          
  			\includegraphics[scale=0.25]{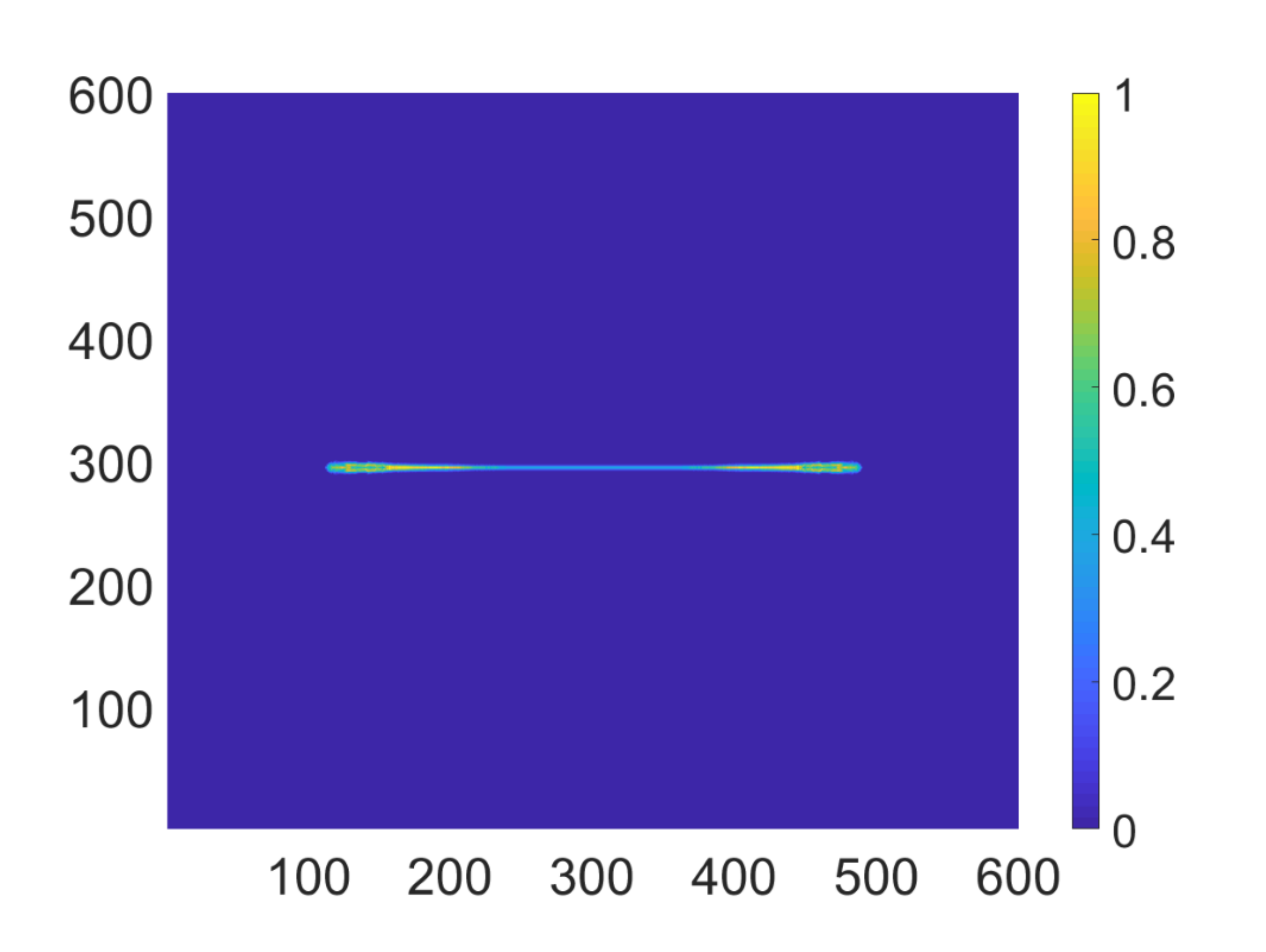}  
  		\end{minipage}
  	}	
  	\subfigure[]
  	{
  		\begin{minipage}{7.5cm}
  			\centering      
  			\includegraphics[scale=0.25]{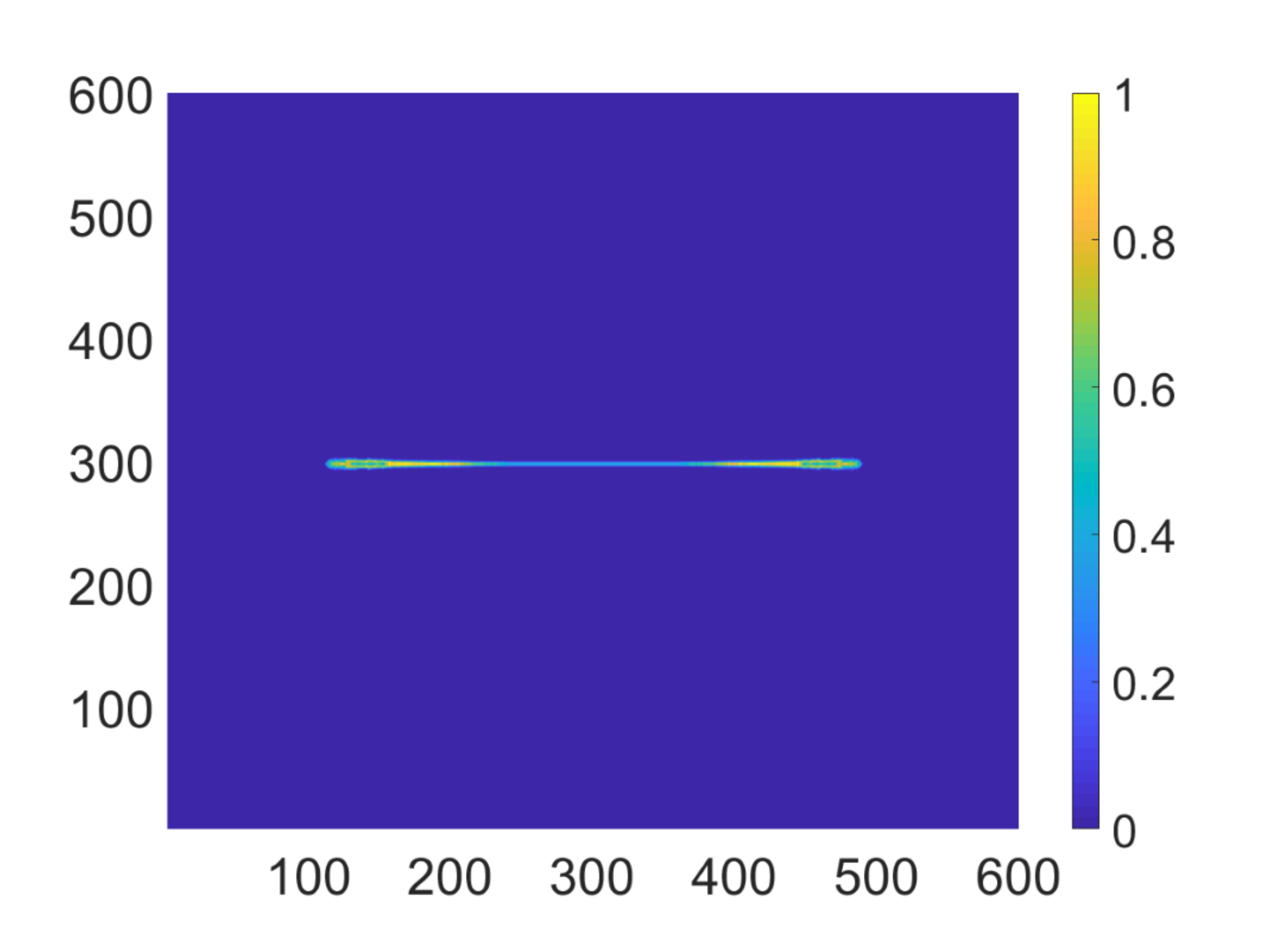}   
  		\end{minipage}
  	}	
  	\caption{Crack simulation results at $1.6708\times10^{-3}s$ with two methods: (a) MSBFM (b)meshfree method} 
  	\label{Fig:9}
  \end{figure}
\begin{table}[!h]
		\caption{Performance of meshfree method and MSBFM with  in 2D model}
		\vspace{0.15in}
		\centering
		\setlength{\tabcolsep}{1mm}{
			\begin{tabular}{ccccccc} \hline
				\multicolumn{2}{c}{Mesh} &$600\times600$&$1200\times1200$&$1800\times1800$&$2400\times2400$&$3000\times3000$\\\hline
				\multicolumn{2}{c}{Time}&1250&1250&1250&1250&1250\\ \hline				
				 \multirow{3}*{Meshfree}
				 &Matrix assembly&21m36s&5h48m&1d1h&3d11h&-\\ 
				 &Time stepping&13m55s&5h50m&1d12h&4d14h&-\\ 				
	
				&Crack factor&8m23s&3h30m&21h4m&2d18h&-\\ \hline
				\multirow{3}*{MSBFM}  & Matrix assembly&2m11s&28m54s&2h28m&8h11m&18h46m\\
				&Time stepping &3m22s&16m20s&1h5ms&4h50m&18h36m\\			
				&Crack factor&8m23s&3h30m&21h4m&2d18h&4d24m\\ \hline
		\end{tabular}}\label{tab:4}
\end{table}

The computational time by using various discretization sizes is shown in Table \ref{tab:4}, and we only calculate the results within 10 days. Here the simulation process is divided into three phases. Phase I: matrix assembly; Phase II:Time stepping;  Phase III: Crack factor. The phase I and the phase II are the same as the non-fracture problem mentioned in Section 5.1, and Phase III mainly includes the calculation of history-dependent scalar-valued $\mu$ and matrix $\mathbf{D}^f$. At each phase, there are differences in the time of non fracture problems and fracture problems.

In Phase I, the advantages of constructing the MSBFM stiffness matrix were retained. Although an additional matrix $\mathbf{D}^f$ is introduced, the elements of matrix $\mathbf{D}^f$ can be obtained from the entries of matrix $\mathbf{K}$, so no additional assembly is required. But we still need to consider the information of material points near the fracture to calculate s when traversing the horizon. Since the crack shape is impossible to estimate, we often need to consider all the material point information in this part, which leads to MSBFM can not save the calculation amount in this part. In actual calculation, we can only consider a preset region and a region with incomplete horizon if all cracks do not exceed this preset region, thus reducing the traversal time.

In Phase II, the matrix $\mathbf{D}^f$ mentioned in (\ref{df}) is computed in each time step, which causes that the computational time of MSBFM in phase I is not strict $O(N\log N)$. $\mathbf{D}^f$ will also affect the time of matrix assembly, but compared with the meshfree method, MSBFM still has computational advantages in the fracture problem because the computational complexity of $\mathbf{D}^f$ does not exceed $O(N^2)$ according to the above analysis.

In the non fracture problem, we do not need to consider Part III. But in fracture problem, the calculation time for solving $\mu$, which is the main part of part III takes up a large part, which is mainly caused by the elongation $s$. $s$ is obtained by a nonlinear form $s=(|\mathbf{x}^{\prime}+\mathbf{u}^{\prime}-\mathbf{x}-\mathbf{u}|-|\mathbf{x}^{'}-\mathbf{x}|)/(|\mathbf{x}^{'}-\mathbf{x}|)$, so it cannot be solved by FFT, which causes that the calculation of this part is usually $O(N^2)$. In the meshfree and MSBFM methods, the calculation amount for solving $s$ is the same.

 \subsection{Peridynamic in 3D body under displacement constraints}
   We perform PD simulations using a 3D model with incomplete horizons and displacement constraints in this section.

  As shown in Fig. \ref{Fig:10}, a block with length $L=1.0\mathrm{~m}$,  width $W=0.3\mathrm{~m}$, and thickness $H=0.3\mathrm{~m}$ is introduced  . Horizon size $\delta$ is chosen as $0.00315\mathrm{~m}$. External loading $b_p=p_0W/h^2$ is applied to the area $\mathcal{D}_s$, and the value of $p_0$ is the same as that in Section 5.1.

  The displacement constraints is imposed on the virtual boundary layer $\mathcal{D}_c$, which means:
     \begin{equation}
   \mathbf{u}_x(x_p,t)=\mathbf{u}_y(x_p,t)=\mathbf{u}_z(x_p,t)=0,\quad \mathbf{x}_p \in A_c
   \end{equation}  
   Elastic modulus $E$ and Poisson's ratio $\nu$ are taken as $2\times 10^5$Mpa and $1/4$, respectively. The material points affected by displacement constraints are treated as material points on $\Omega_c$.
     \begin{figure}[h]
   	\centering    	
   	\subfigure[]
   	{
   		\begin{minipage}{6cm}
   			\centering          
   			\includegraphics[scale=0.5]{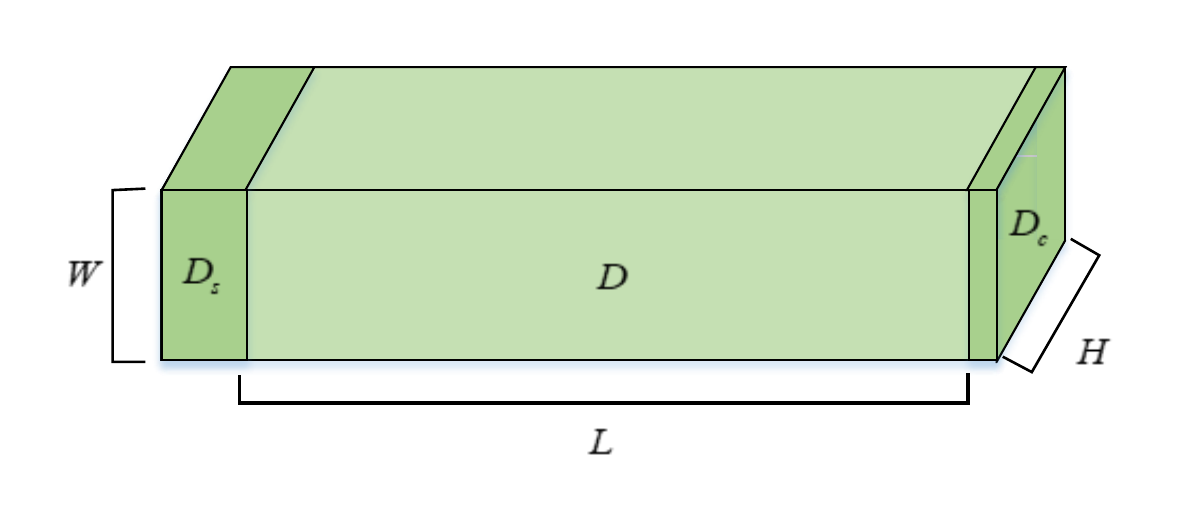}  
   		\end{minipage}
   	}	
   	\subfigure[]
   	{
   		\begin{minipage}{6cm}
   			\centering      
   			\includegraphics[scale=0.5]{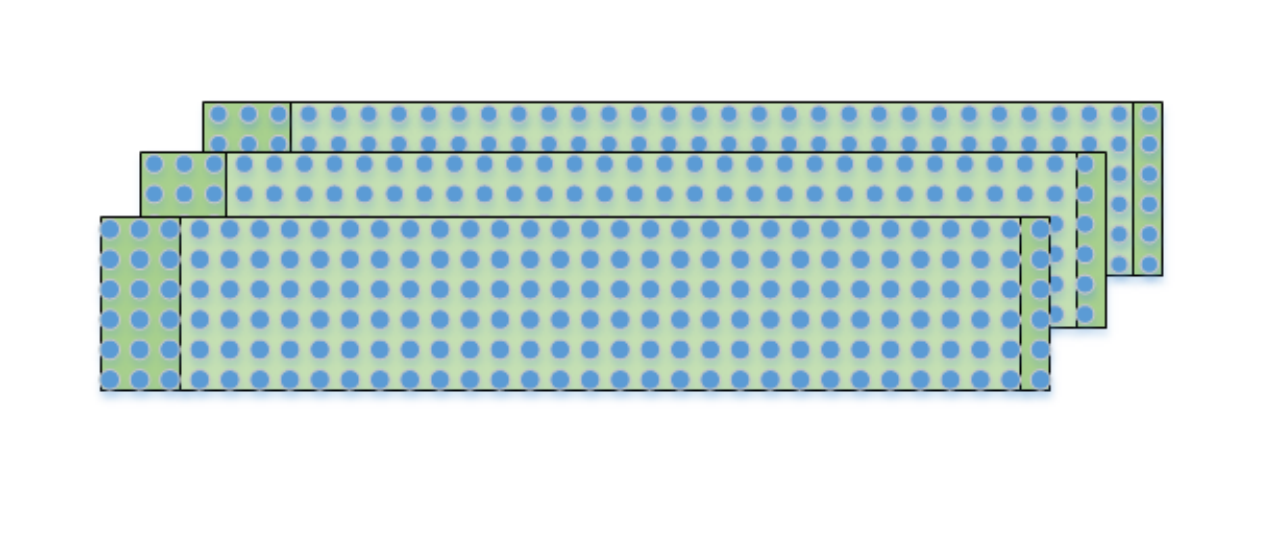}   
   		\end{minipage}
   	}	
   	\caption{Geometry of a block under displacement constraints and its discretization: (a)block with external loading and displacement constraints; (b)a simple example of uniform grid structure.} 
   	\label{Fig:10}
   \end{figure}

The grid of size $100\times30\times30$ is used to discretize the model and perform time stepping through the ADR method. The displacement variations of the two algorithms are provided in in Fig.  \ref{Fig:11}, which verify the accuracy of MSBFM in 3D problems.
          \begin{figure}[h]
  	\centering    	
  	\subfigure[]
  	{
  		\begin{minipage}{6cm}
  			\centering          
  			\includegraphics[scale=0.25]{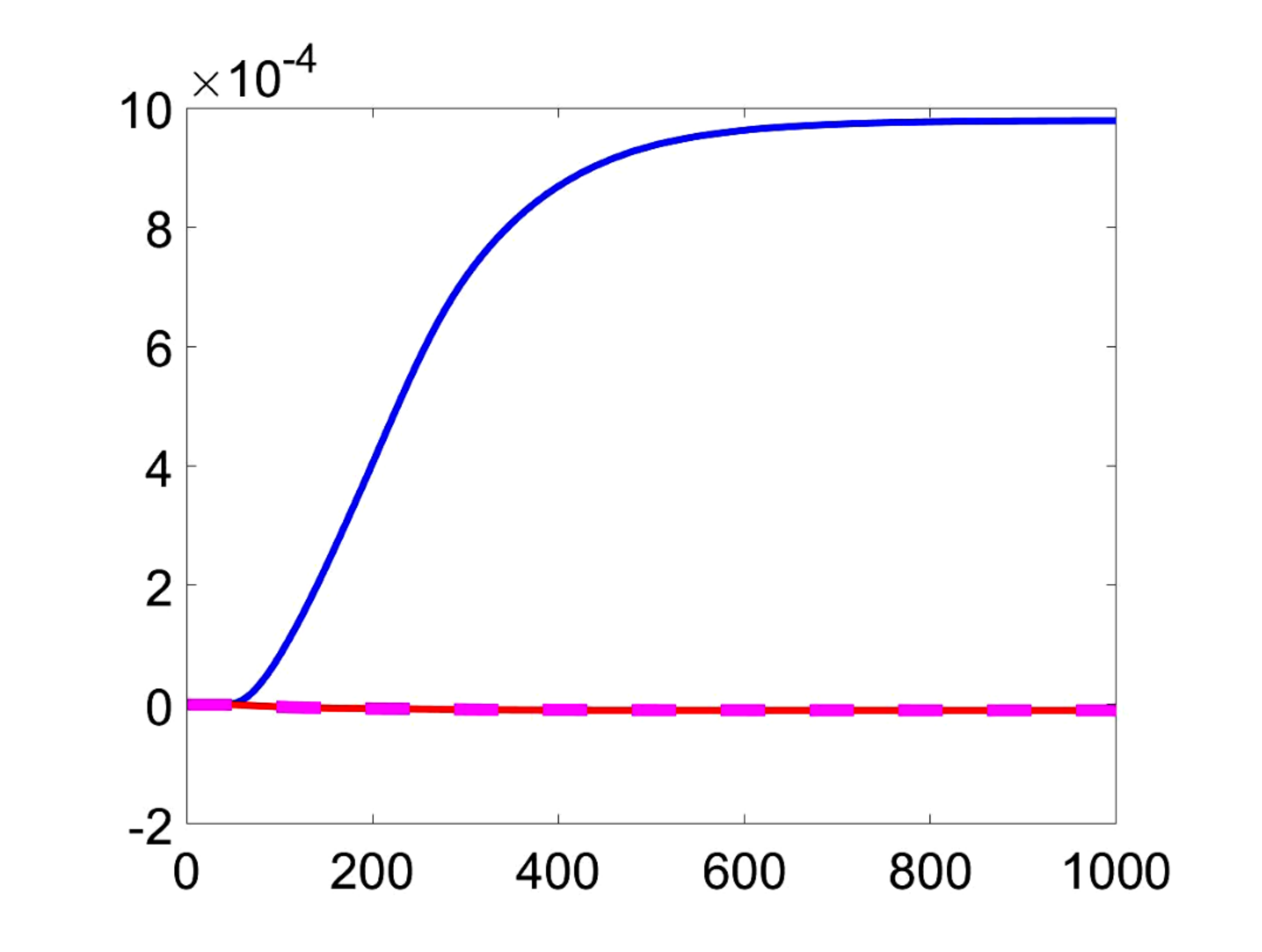}  
  		\end{minipage}
  	}	
  	\subfigure[]
  	{
  		\begin{minipage}{6cm}
  			\centering      
  			\includegraphics[scale=0.25]{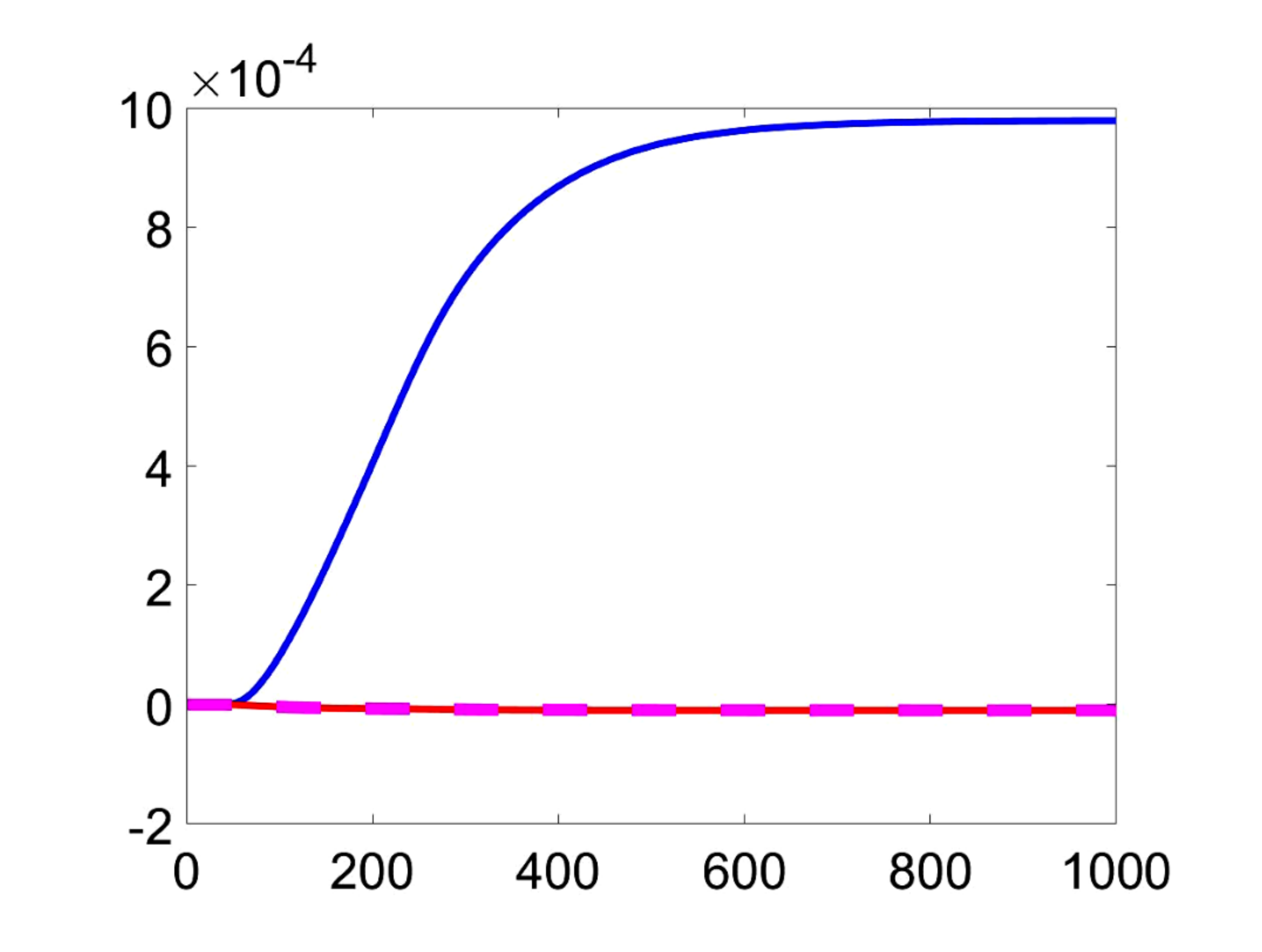}   
  		\end{minipage}
  	}	
  	\caption{Displacement variation with respect to step number in $x$ direction (blue), $y$ direction(red), $z$ direction(green): (a)MSBFM (b)meshfree method} 
  	\label{Fig:11}
  \end{figure}

Table \ref{tab:6} compares the computational time required to perform the simulations using MSBFM and the meshfree PD method, and Only results not exceeding 10 days are considered. Similar to the two-dimensional problem, the computational time of the three-dimensional non fracture problem can also be divided into two phases: phase I: Compute $\mathbf{f}$, including the calculation of $\mathbf{f}^x$, $\mathbf{f}^y$, $\mathbf{f}^z$ and the form $\mathbf{Du}$; B: Matrix assembly, including searching the horizon of the material point and initialization of matrices $\mathbf{A}^{xx}$, $\mathbf{A}^{xy}$, $\mathbf{A}^{xz}$, $\mathbf{A}^{yx}$, $\mathbf{A}^{yy}$, $\mathbf{A}^{yz}$, $\mathbf{A}^{zx}$, $\mathbf{A}^{zy}$, $\mathbf{A}^{zz}$ and corresponding matrices $\mathbf{D}^e$.
 \begin{table}\begin{center}
 		\caption{Performance of meshfree method and MSBFM  in 3D model}
 				\vspace{0.15in}
 		\setlength{\tabcolsep}{3mm}{
 			\begin{tabular}{ccccc} \hline
 				\multicolumn{2}{c}{Mesh}&$100\times30\times30$&$200\times60\times60$&$300\times90\times90$\\ \hline
 				\multicolumn{2}{c}{Time}&1000&2000&3000\\\hline
 				
 				\multirow{2}*{Meshfree}	& Matrix assembly&1m40s&1h47m&4d17h\\ &Time stepping&35m5s&1d13h&-\\ 
 			    \hline
 				\multirow{2}*{MSBFM }&Matrix assembly&37s&20m12s&3h20m\\&Time stepping &3m2s&49m15s&9h10m\\
 				 \hline
 		\end{tabular}}\label{tab:6}
 \end{center}\end{table}

For the part I, the calculation amount of $\mathbf{D}^e$ can be obtained from $N^2(4\delta/h)=O(N^{\frac{7}{3}})$. In the part of traversing horizon of material points, the ratio of the computational time is the same for the 3D model and the 2D model. In fact, the the computational time of part I  depends on the number of material points rather than the dimension of the point, and in three-dimensional problems, the number of material points tends to be more than in two-dimensional problems, thus MSBFM will be more advantageous in  part I in three-dimensional problems.

 We observed a high rate between MSBFM and meshfree method in part II, and this is the result of the differences in dimensions between the 2D model and the 3D model. When the $N_x$, $N_y$, $N_z$ increases to $\beta N_x$, $\beta N_y$, $\beta N_z$, $N$ increases to $\beta^3 N$. Similar to the analysis in section 5.1, the computational time of our algorithm will increase by $\beta^3$ times, while the meshfree method will increase by $\beta^6$ times. This is why the time ratio of part II is larger than that of 2D model.

One problem that needs to be noted is that the proportion of material points on $\Omega_e$ in the total material points in the 3D model will also increase. Therefore, when the surface correction algorithm is adopted, the computational advantage of algorithms may not be obvious.

\subsection{Kalthoff-Winkler experiment} 
To simulate the calculation rate of MSBFM algorithm on 3D fracture model, a KW example is introduced in this subsections, as shown in Fig. (\ref{Fig:12}).

The problem description is as follows: A block of length $L=0.2 \mathrm{~m}$, width  $W=0.1 \mathrm{~m}$, and thickness  $h=0.009\mathrm{~m}$ with two thin notches is subjected to the incomplete horizons and impactor. The notch in this model has width $h_0=0.0015 \mathrm{~m}$, length $a_0=0.05 \mathrm{~m}$, and distance between notches $d$ are $0.05 \mathrm{~m}$ and the impactor's diameter $D$ and height $H$ are all $0.05 \mathrm{~m}$. Horizon size $\delta$, density $\rho$, elastic modulus $E$ and Poisson's  rate $\nu$ considered are chosen as in Section $5.3$. Initial velocity $v_0=32 \mathrm{~m} / \mathrm{s}$ is imposed to ensure crack propagation. Critical stretch $s_0$ is chosen for judging fracture, which is taken as $0.01$. We conduct this simulations by using Velocity Verlet, and each time step is $1.3367\times10^{-8}$.

\begin{figure}[h]
	\centering    	
	\subfigure[]
	{
		\begin{minipage}{6cm}
			\centering          
			\includegraphics[scale=0.3]{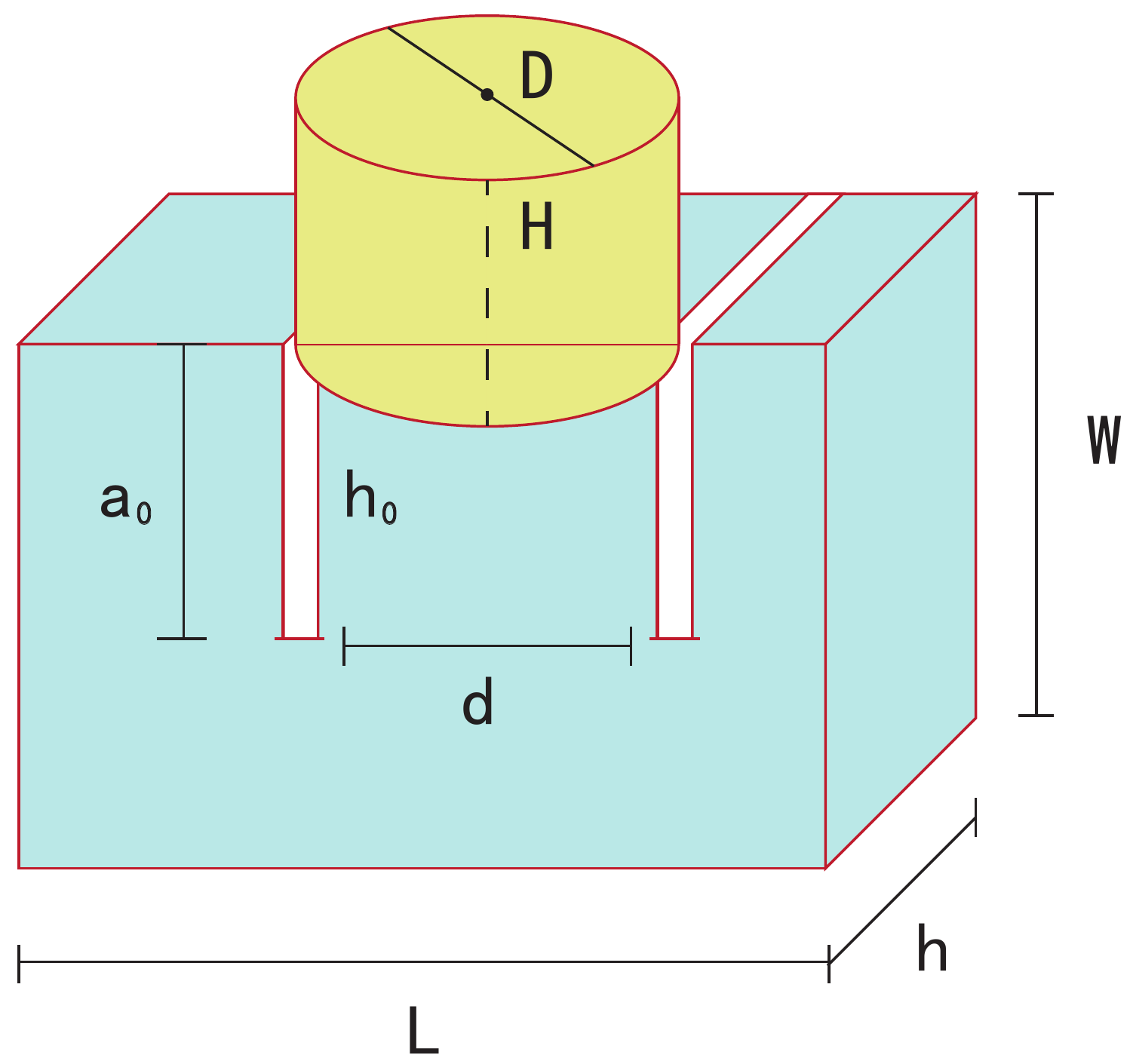}  
		\end{minipage}
	}	
	\subfigure[]
	{
		\begin{minipage}{6cm}
			\centering      
			\includegraphics[scale=0.3]{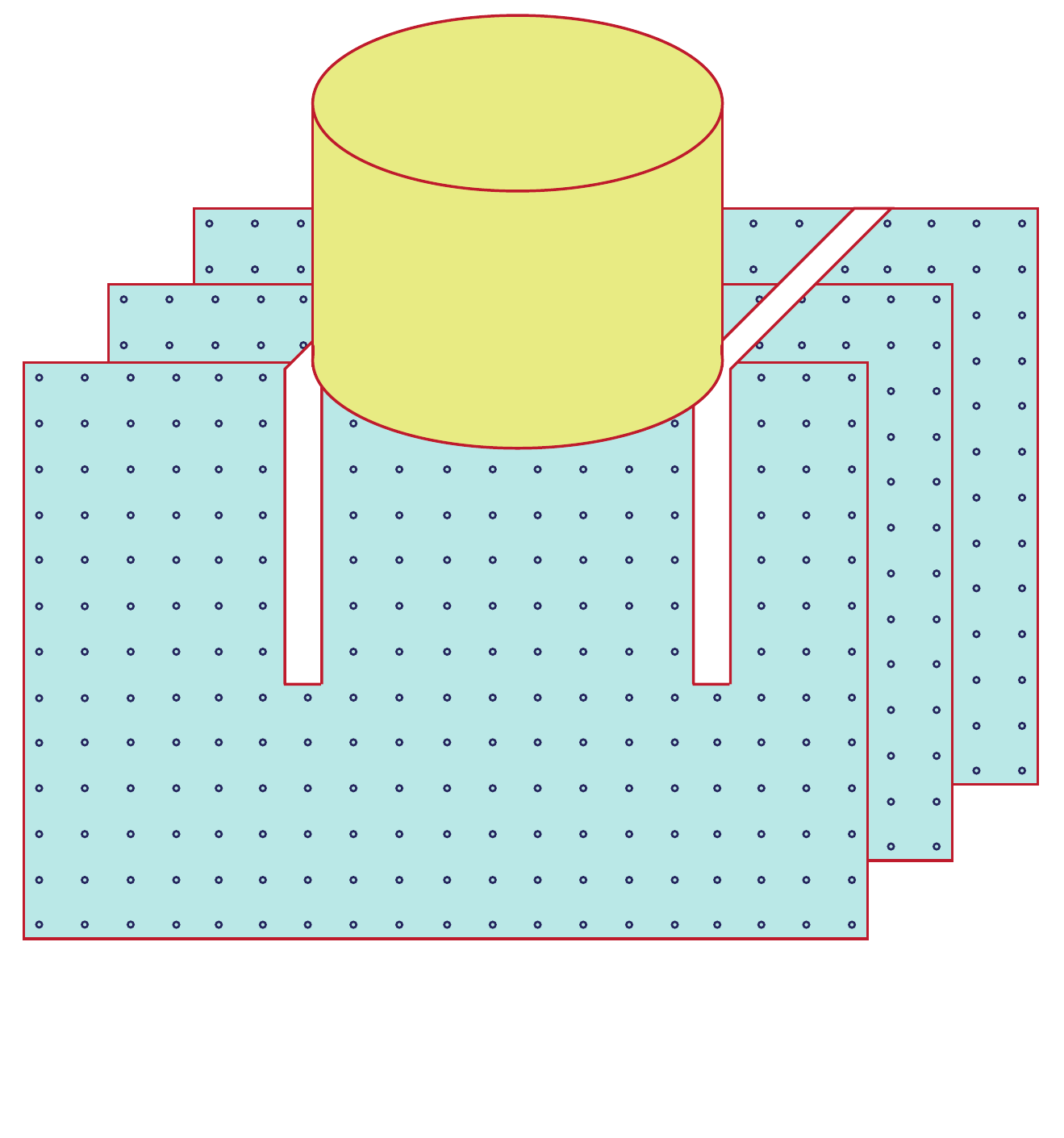}   
		\end{minipage}
	}	
	\caption{Geometric of Kalthoff-Winkler test and its discretization(a)block used in  Kalthoff-Winkler test; (b)a simple example of uniform grid structure.} 
	\label{Fig:12}
\end{figure}

The crack propagations by our algorithm and meshfree  method are provided in Fig. \ref{Fig:13}. we can obverse that the cracking angle is $135^{\circ}$, which is is consistent with the results obtained in the standard experiment. The crack propagations in our algorithm at different time steps shown in Fig. \ref{Fig:14} and the computational time is shown in Table \ref{tab:7}. For more complex meshes, both FMBM and meshfree methods have exceeded the time limit, which is mainly caused by point arrangement. If other point arrangement methods are used, the calculation time will be greatly reduced.
\begin{figure}[h]
	\centering   	
	\subfigure[]
	{
		\begin{minipage}{7cm}
			\centering          
			\includegraphics[scale=0.2]{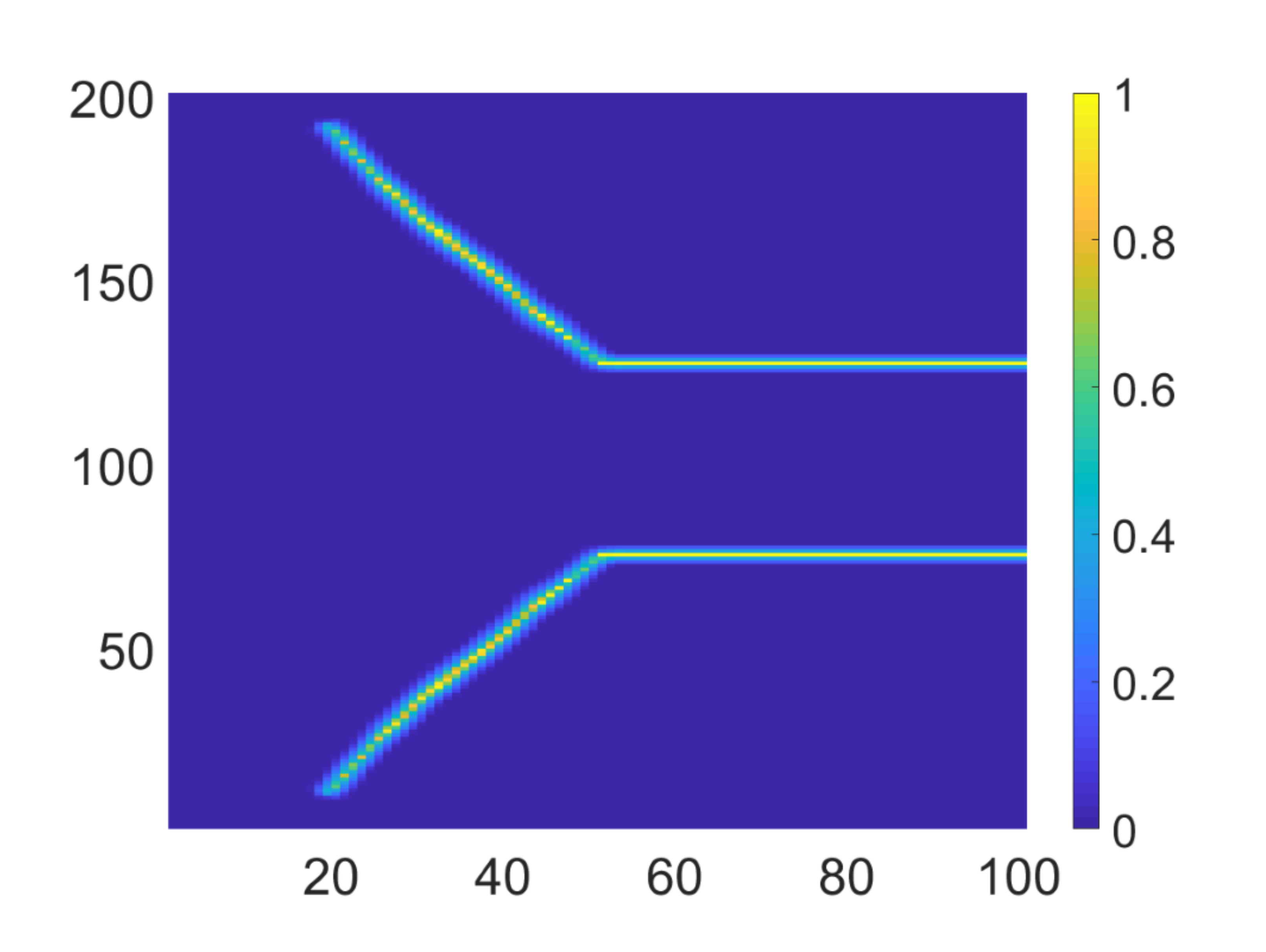}  
		\end{minipage}
	}	
	\quad
	\subfigure[]
	{
		\begin{minipage}{7cm}
			\centering      
			\includegraphics[scale=0.2]{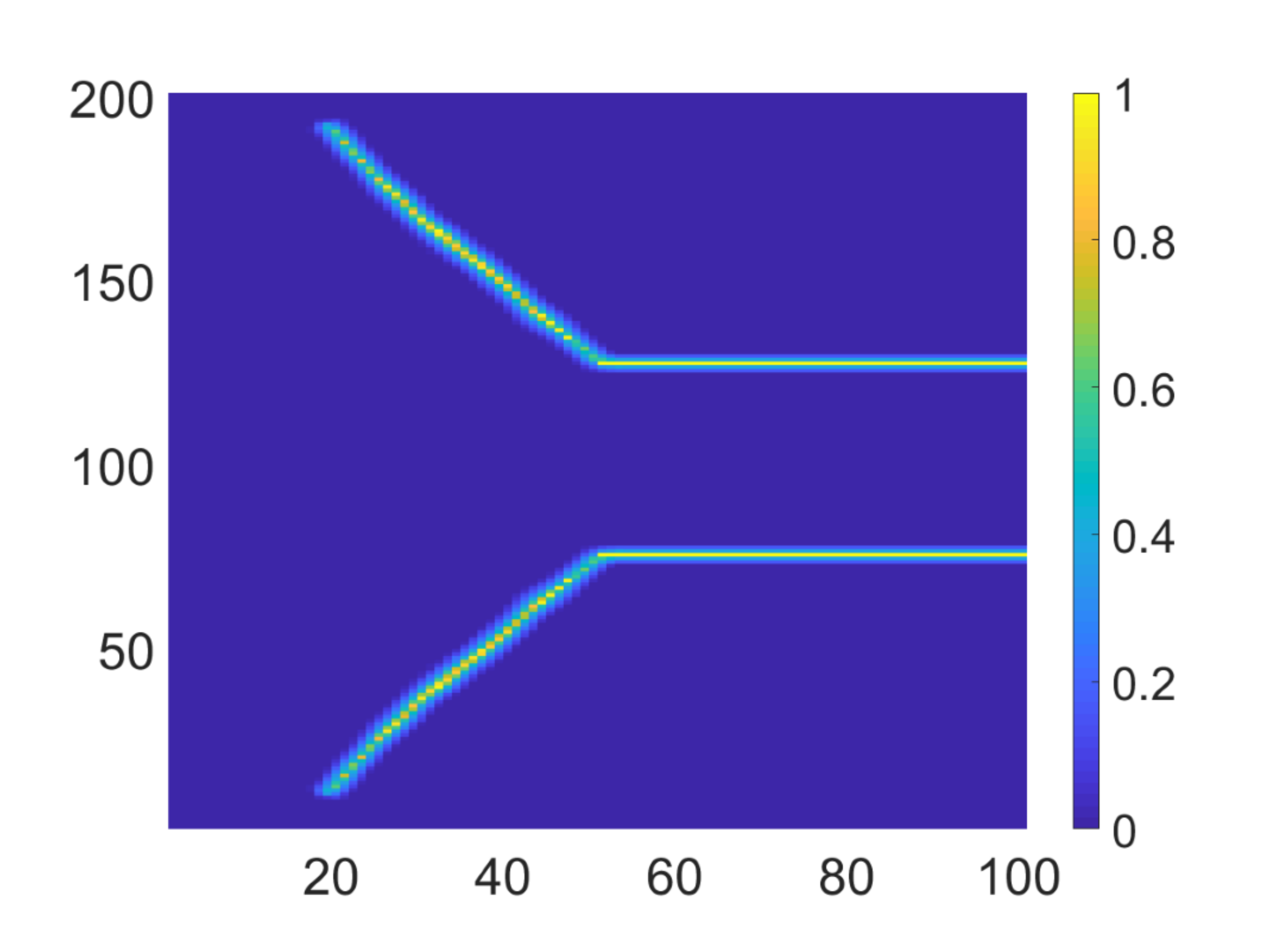}   
		\end{minipage}
	}	
	\caption{crack simulation results at $1.1754\times10^{-8}s$ with two methods: (a)MSBFM (b)meshfree method} 
	\label{Fig:13}
\end{figure}
\begin{figure}[!htb]
	\centering   	
	\subfigure[]
	{
		\begin{minipage}{.3\linewidth}
			\centering          
			\includegraphics[scale=0.15]{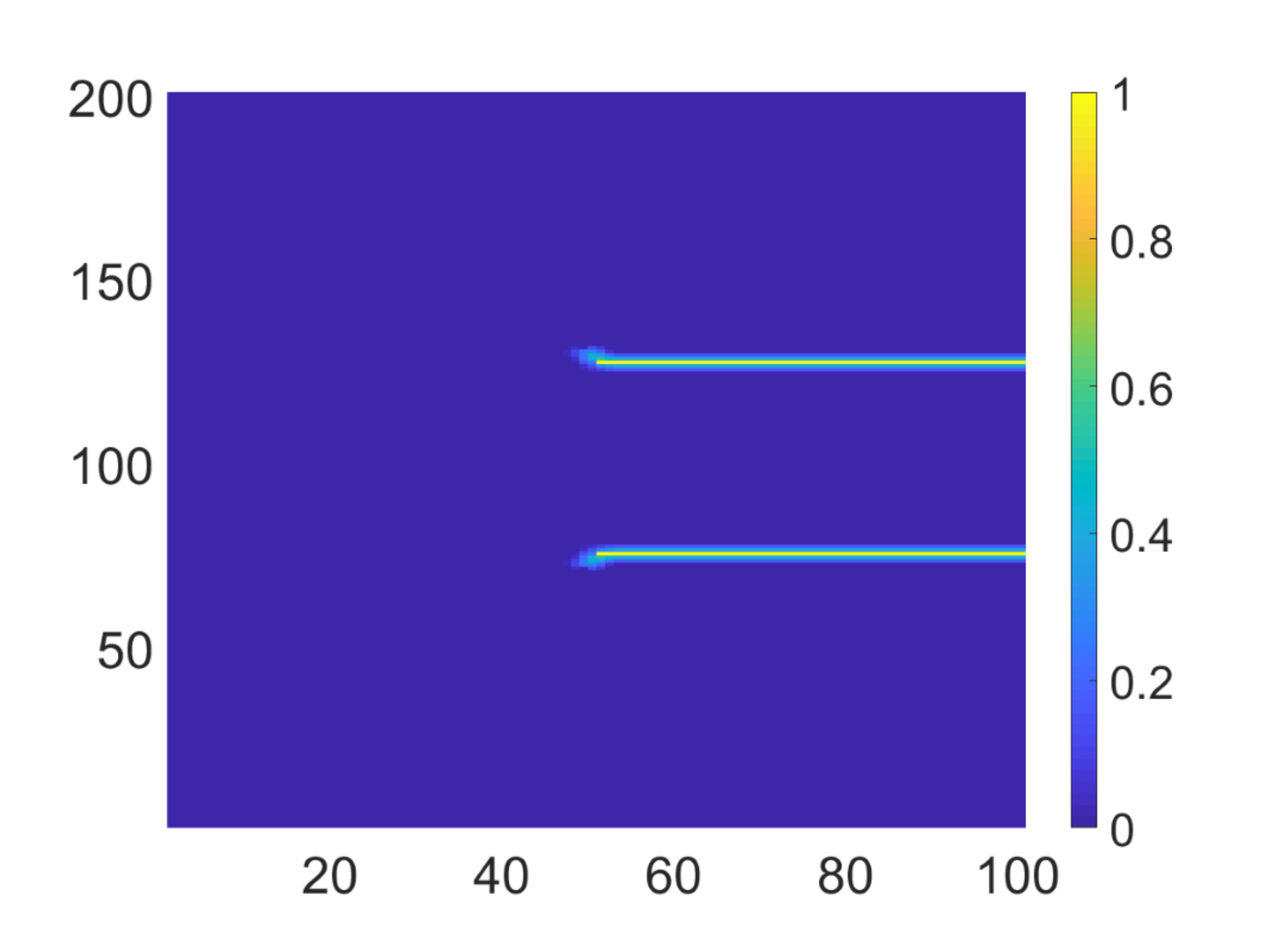}  
		\end{minipage}
	}	
	\subfigure[]
	{
		\begin{minipage}{.3\linewidth}
			\centering      
			\includegraphics[scale=0.15]{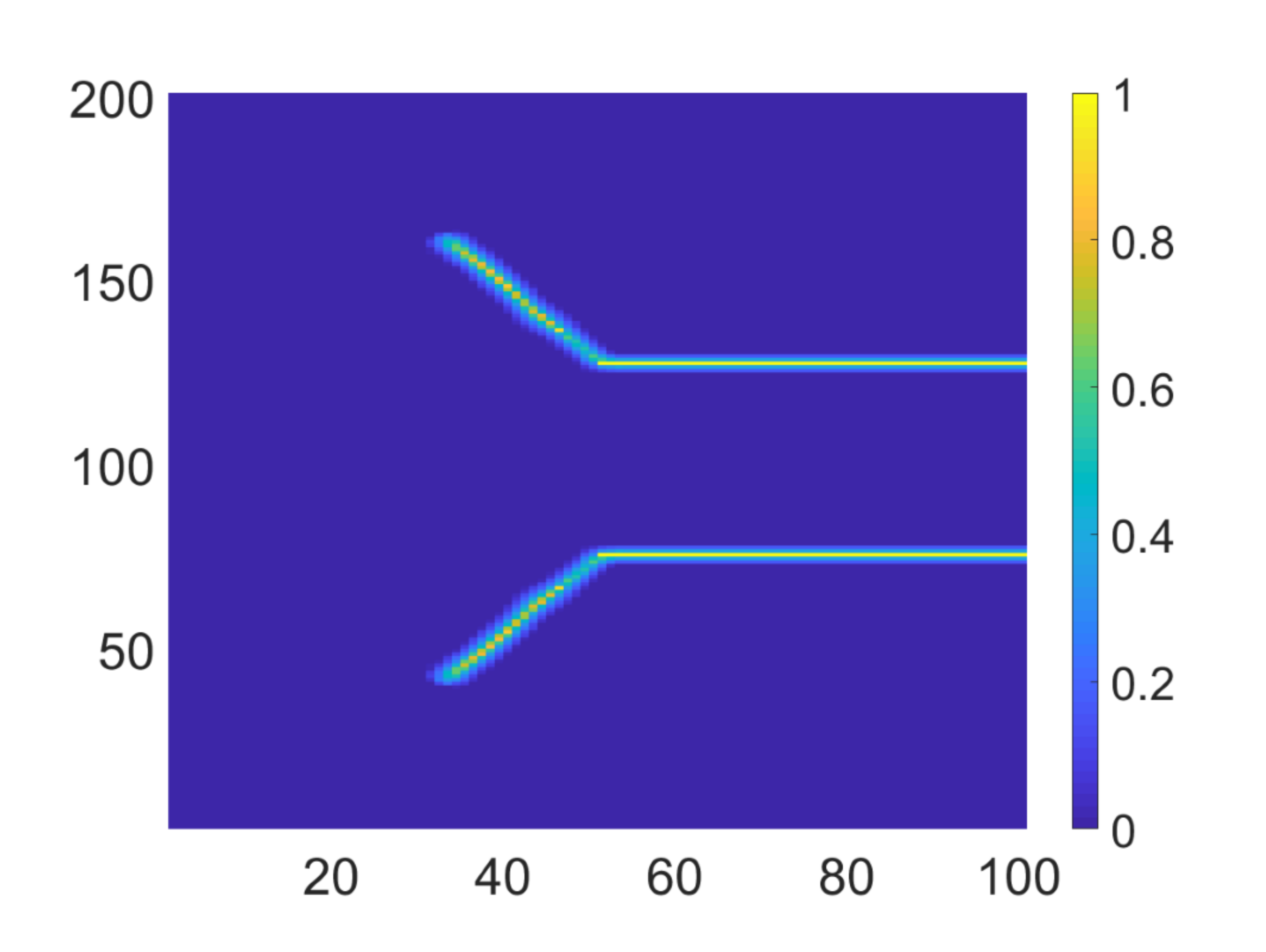}   
		\end{minipage}
	}
	\subfigure[]
	{
		\begin{minipage}{.3\linewidth}
			\centering      
			\includegraphics[scale=0.15]{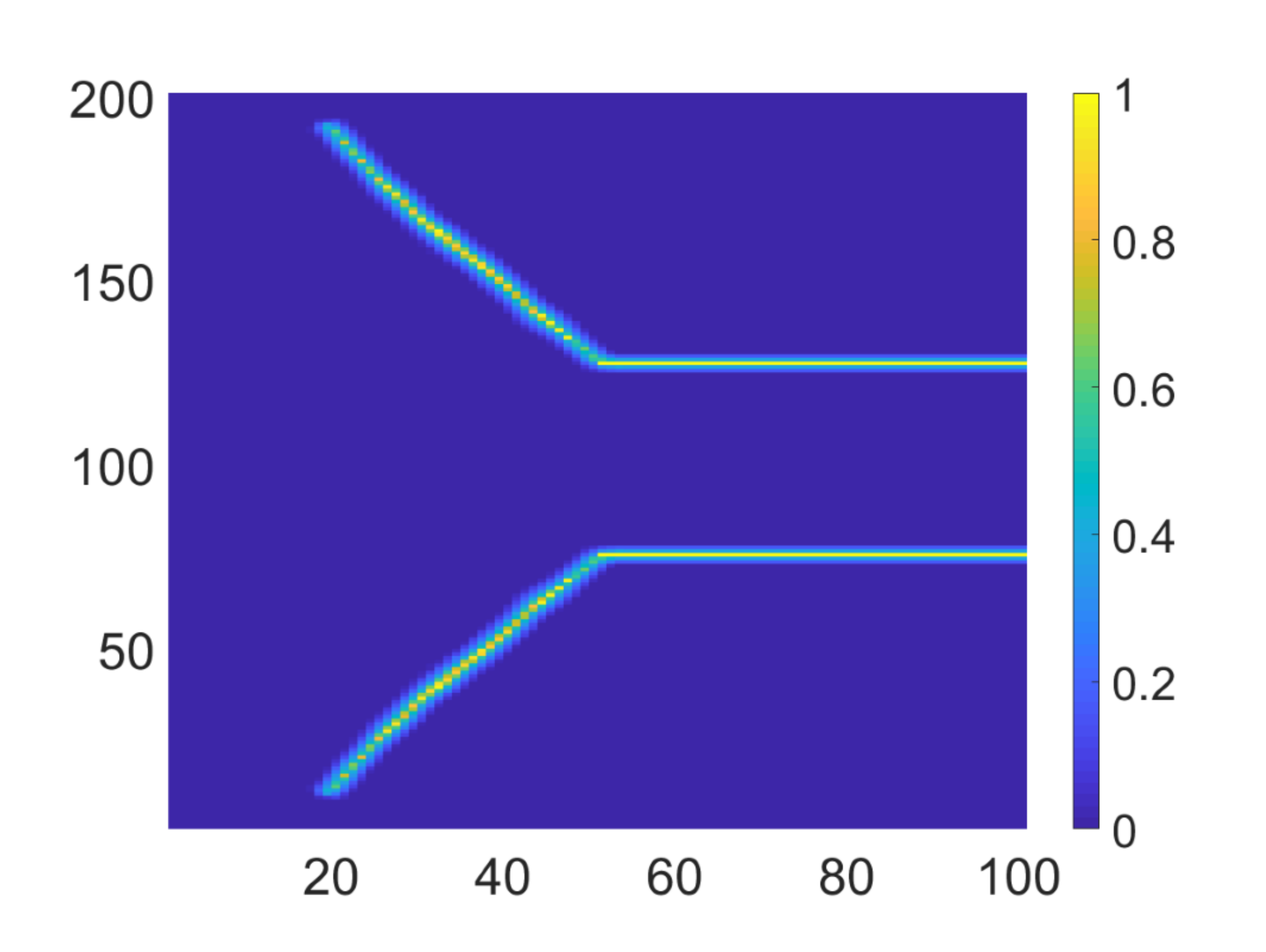}   
		\end{minipage}
	}	
	\caption{crack simulation results by MSBFM (a)at 350 time steps(b)at 650 time steps(c)at 1350 time steps} 
	\label{Fig:14}
\end{figure}
\begin{table}\begin{center}
		\caption{Performance of meshfree method and MSBFM  in KW experiment}
				\vspace{0.15in}
		\setlength{\tabcolsep}{4mm}{
			\begin{tabular}{cccccc} \hline
				\multicolumn{2}{c}{Mesh}&$201\times101\times9$&$402\times202\times18$\\\hline
				\multicolumn{2}{c}{Time}&1350&1350\\ \hline				
				\multirow{3}*{Meshfree}&Matrix assembly&9m25s&11h10m&& \\&Time stepping&1h16m&1d4h&\\ 
				&Crack factor&40m30s&17h5m&&\\ \hline
				\multirow{3}*{MSBFM} 	&Matrix assembly &6m29s&7h39m&&\\&Time stepping&13m27s&2h5m&\\
				&Crack factor&40m25s&17h5m&&\\ \hline
		\end{tabular}}\label{tab:7}
\end{center}\end{table}

As expected, the KW experiment retains the computational advantage of MSBFM in the process of computing $\mathbf{f}$. However, the number of broken bonds in 3D model will also increase, which mainly affects two aspects: first, the construction of matrix $\mathbf{D^f}$, whose calculation is usually $O(N^2)$, which will affect the process of computing $\mathbf{f}$ and matrix assembly. The other is the calculation of $\mu$, because of the increase of broken bonds, the calculation of fracture simulation will also increase, which occupies the main part in each time step. 

In fact, since the number of broken bonds cannot exceed the total number of bonds, the computational advantages of MSBFM can be maintained in most fracture models, especially in 3D models. The fewer broken bonds, the more obvious the computational advantage of the MSBFM algorithm.

Comparisons between the computational efficiency of the new MSBFM method for PD models with that of the original meshfree discretization of PD formulations by four examples showed the computational and storage advantages of our algorithm, especially in 3D problems. One
can now easily simulate fracture problems by selecting some material points instead of all material points by using MSBFM, which reduces memory allocation and maintains high accuracy compared with the mesh free method

\section{Conclusions}
In this paper, we introduce a matrix-structure-based fast method(MSBFM). In this framework, the stiff matrix is decomposed into a summation of several matrices according to the model's other boundary conditions and fracture conditions. Following these decompositions, FFT and its inverse operation are used to calculate the PD integral with the cost of $O(N\log N)$ instead of $O(N^2)$ required by the usual meshless or FEM discretization methods of the PD model. Because of the Fourier transform, storing all the information about the material points and their horizon is no longer necessary, thus reducing the storage cost from $O(N^2)$ to $O(N)$ of the meshfree or FEM discrete method. Therefore, the time for initializing the matrix is also reduced. For the time-dependent problems and quasi-static problems, the time marching schemes are used to simulate.

The method mentioned in this paper applies to most nonlocal models as long as their discrete forms can be written in a matrix. This paper focuses on the bond-based PD model, and the following numerical test are performed: two-dimensional non-fracture problems with loading
two-dimensional fracture problems with displacement constraints, three-dimensional non-fracture problems with two kinds of boundary conditions. The results are in good agreement with the theory. The comparison with the computational speed of the meshfree method shows that MSBFM can reduce the computational time of tens of days in this method to several hours. This means that for complex fracture problems, the selection of PD nodes and the computational cost is no longer the main obstacles for complex fracture problems.

The algorithm still depends on the matrix structure to some extent, which means it must be a quasi-Toeplitz structure. In other words, most entries satisfy the Toeplitz structure. Efforts are underway to extend the application of the SFPD algorithm to more PD problems with a complex matrix structure, including state-based problems, coupling problems, and nonlinear problems.

\section*{Acknowledgements}
The work was carried out at Marine Big Data Center of Institute for Advanced Ocean Study of Ocean University of China.

\end{document}